\documentclass[a4paper,11pt]{article}

\usepackage{graphicx,float}
\usepackage{amsfonts,latexsym,amssymb} 
\usepackage{mathrsfs,amsmath}
\usepackage{color}
\usepackage{soul}
\newcommand{\vecg}{\mathbf{g}}

\newcommand{\vecu}{\mathbf{u}}
\newcommand{\vecv}{\mathbf{v}}
\newcommand{\vecw}{\mathbf{w}}

\newcommand{\vecz}{\mathbf{z}}

\newcommand{\Proof}{{\noindent {\bf Proof:   }}}
\newcommand{\EndProof}{{\hfill $\Box \qquad $ \endtrivlist}\par \vspace{0.5cm}}%

    \newtheorem{theorem}{Theorem}[section]
    \newtheorem{lemma}[theorem]{Lemma}
    \newtheorem{proposition}[theorem]{Proposition}
    \newtheorem{definition}{Definition}[section]
   \newtheorem{hypo}{Assumption}[section]
    \newtheorem{remark}{Remark}[section]

\newcommand{\CorrR}[1]{{\color{black} #1}}
\newcommand{\CorrB}[1]{{\color{black}#1}}

\begin{document}
\title{Hysteresis and controllability of affine driftless systems: some case studies\footnote{The work has been partially developed within the OptHySYS project of the University of Trento and also partially financially supported by  Gruppo Nazionale Analisi Matematica Probabilit\`a e Applicazioni (GNAMPA).}}
\date{ }
%
\maketitle

\begin{center}
Fabio Bagagiolo, \footnote{Dipartimento di Matematica, Universit\`a degli studi di Trento 
              via Sommarive,14 38123 Povo (TN) Italy, fabio.bagagiolo@unitn.it}
Marta Zoppello\footnote{Politecnico di Torino, Corso Duca degli Abruzzi, 24
10129 Torino, Italy, marta.zoppello@polito.it}\end{center}
%
%
\begin{abstract} We investigate the controllability of some kinds of driftless affine systems where hysteresis effects are taken into account, both in the realization of the control and in the state evolution. In particular we consider two cases: the one when hysteresis is represented by the so-called play operator, and the one when it is represented by a so-called delayed relay. In the first case we prove that, under some hypotheses, whenever the corresponding non-hysteretic system is controllable, then we can also, at least approximately, control the hysteretic one. This is obtained by some suitably constructed approximations for the inputs in the hysteresis operator. In the second case we prove controllability for a generic hysteretic delayed switching system. Finally, we investigate some possible connections between the two cases.\end{abstract}
%
%

\section*{Introduction}
\label{intro}
The study of different kinds of mechanical systems provides a rich area for mathematical investigation, and vice-versa mathematics may enlighten mechanical phenomena. In particular, we may have possible implications on the design of artificial devices that could be used in different context, from medicine to industry. Nowadays relying on technology has become fundamental and the mathematical modelization of mechanics underlying any real system is crucial for the development of any sophisticated technology. In particular, such models ought to contain a control, so that control theory is likely the appropriate mathematical framework for this issue. Furthermore, many of the models used to describe real mechanical (as well as physical, biological, economic and social) systems may present an intrinsic memory phenomenon.  To take into account this particular memory behavior, one may introduce
a suitable memory-term and pursue controllability results in this enlarged setting.
Often, such a memory effect is of the so-called rate-independent type. This means that the actual state of the system depends on its whole past history via the sequence of reached states only, independently of the time-scale.  Sometimes this behavior may be seen as a sort of delay in the reaction to some external forces (as well as to some external controls). This phenomenon is known as hysteresis, and just to name a few of examples, besides the classical ferromagnetic theory, we quote hysteresis in phase transitions (see Brokate-Sprekels \cite{brospr}), hysteresis in filtration through porous media (Bagagiolo-Visintin \cite{bagvis}), hysteresis in economics (see Gocke \cite{goc}), hysteresis in transmission and consensus problems (see Ceragioli et al. \cite{CPF}). The mathematical studies of hysteresis phenomena as functional operators, representing the input-output hysteresis relationship, were introduced by Krasnosel'ski\v{\i} and his co-worker \cite{KR} (see also Visintin \cite{Visintin94}). 
This kind of operators are non linear and non differentiable, even if some possible definitions of derivatives were given, see for example Brokate- Krej\v{c}\'i \cite{brokre}. Anyway, those definitions essentially involve derivatives of the output with respect to time. When spatial derivatives must be taken into account, 
the dependence on the past history is, up to the knowledge of the authors, an unsolved problem. This fact, in the controllability setting,  prevents the use  of local techniques, for which the application of classical tools in geometric control theory, as Chow theorem, Lie brackets and so on (Chow \cite{Chow}, Coron \cite{Coron56}, Bressan \cite{Bressan07}, Agrachev et al. \cite{Agrachev12}) is not immediate. Indeed such tools involve the spatial derivatives of the dynamic vector fields which, as already said, seem to be meaningless in the presence of hysteresis. Specific studies are then required but they are not well presented in the literature, despite the importance of the problem.

Due to the difficulties described above, in this paper we assume the controllability of the system without hysteresis and study the case when that system is perturbed by a hysteresis effect. This kind of situation is also common in the applications, for example in the case of switching systems (see Liberzon \cite{lib}, see also the recent Bagagiolo at al. \cite{BagagioloZoppello16}) or in the case of stabilization of systems representing a single input single output plant interconnected with a hysteresis disturbance (see for example Cocetti et al. \cite{coczacbagber} and the references therein). Moreover, the present study seems to be the first attempt in the direction of the control of driftless affine ordinary systems with hysteresis (for a first attempt in the framework of semilinear parabolic partial differential equations, see Bagagiolo \cite{baghystPDE}. See also Gavioli- Krej\v{c}\'i \cite{gavkre}.)

More precisely we focus on the driftless control-affine system in $\mathbb{R}^n$

\begin{equation}
\label{eq:intro_no_hyst}
\begin{cases}
\dot{\vecz}=\sum_{i=1}^{m}\vecg_i(\vecz)u_i\\
\vecz(0)=\vecz_0
\end{cases}
\end{equation}

\noindent
and we assume that the hysteresis effect is described by the so-called  \textit{play operator} $\mathcal P$  (see Visintin \cite{Visintin94}) which maps a continuous time-dependent function $\zeta$ (the input) to a continuous time-dependent functions ${\cal P}[\zeta]$ (the output), or by the delayed relay operator $h_\rho$, whose output is instead a piecewise-constant time-dependent function and which may be used to model situations of discontinuous switching dynamics. These two remarkable examples of hysteresis operators are introduced in Section \ref{sec:operators}, and, as we are going to describe in Subsection 4.4, they are also intimately related.

Regarding the play operator, it can be introduced in system (\ref{eq:intro_no_hyst}) in two different ways:

\begin{equation}
\label{eq:intro}
\begin{cases}
\dot{\vecz}=\sum_{i=1}^{m}\vecg_i(\vecz){\mathcal P}[u_i]\\
\vecz(0)=\vecz_0
\end{cases}
\ \ \ \ \ \ \ 
\begin{cases}
\dot{\vecz}=\sum_{i=1}^{m}\vecg_i({\mathcal P}[\vecz])u_i\\
\vecz(0)=\vecz_0
\end{cases}
\end{equation}

\noindent
On one hand we apply it in the controls (see \eqref{eq:intro}-left), on the other hand we introduce the hysteresis in the state variables \eqref{eq:intro}-right. These two cases may model respectively the situation where the control is performed by an external magnetic field (see for example Alouges at al. \cite{AlougesDeSimone14, AlougesDeSimone17}) and where there could be a sort of lack of information in the state-variable, for example in the synthesis of feedback controls (see for example Bauso et al. \cite{Bauso}, Cocetti et al. \cite{coczacbagber}, Logemann et al. \cite{Logeman} and Tarbouriech et al. \cite{Tarbouriech} for the case of linear systems).

The first case is addressed in Section \ref{sec:hysteresis_in_control}. There, suitably using the properties of the play operator, we obtain an approximate controllability result via the construction of a suitable sequence of continuous controls $u^k$ such that ${\mathcal P}[u^k]$ converges in $L^1$ to the (possible discontinuous) control $\overline u$, good for the non-hysteretic case. 

The second case is addressed in Section \ref{sec:hysteresis_in_state}. We restrict ourselves to a suitable class of triangular systems for which we still construct an approximating sequence of controls generating trajectories converging to the good trajectory for the non-hysteretic case. 

Such a class of systems has strong connection with the so-called Heisenberg systems and Carnot groups of step 2, and we also give a possible mechanical justification for it.

In Section \ref{sec:hysteresis_in_state}, for the case of hysteresis in the state variables, we also analyze  the situation in which the system switches between different dynamics and we model it by the introduction of a delayed relay hysteresis operator.  Such a situation occurs, for example, when there is a change in the dynamics depending on the state (e.g. when crossing some hyperplanes of $\mathbb{R}^n$). It has been successfully introduced to solve chattering problems (see Ceragioli et al. \cite{CPF}) or to get complete controllability results (see Bagagiolo et al. \cite{BagagioloZoppello16}).
We prove controllability for such a problem. 

Finally, in Subsection \ref{S_P_Hysteresis}, we give a controllability result for the case where the hysteresis/memory effect is given by the sum of a finite number of delayed relays. As we will see, this situation is also extremely related to the approximation of the continuous play operator, and hence promising in order to obtain more general controllability results.

\section{Hysteresis operators}
\label{sec:operators}
Hysteresis phenomena often occur in mechanical systems such as gear systems, hydraulic controlled valves or systems governed by an external magnetic field. These systems experience a particular memory effect, the rate independent one which is persistent and scale invariant. In this section we describe the mathematical properties of two operators used to model the hysteresis phenomena. Let us start with the so called play operator.
\begin{figure}[H]
\begin{center}
\includegraphics[scale=0.35]{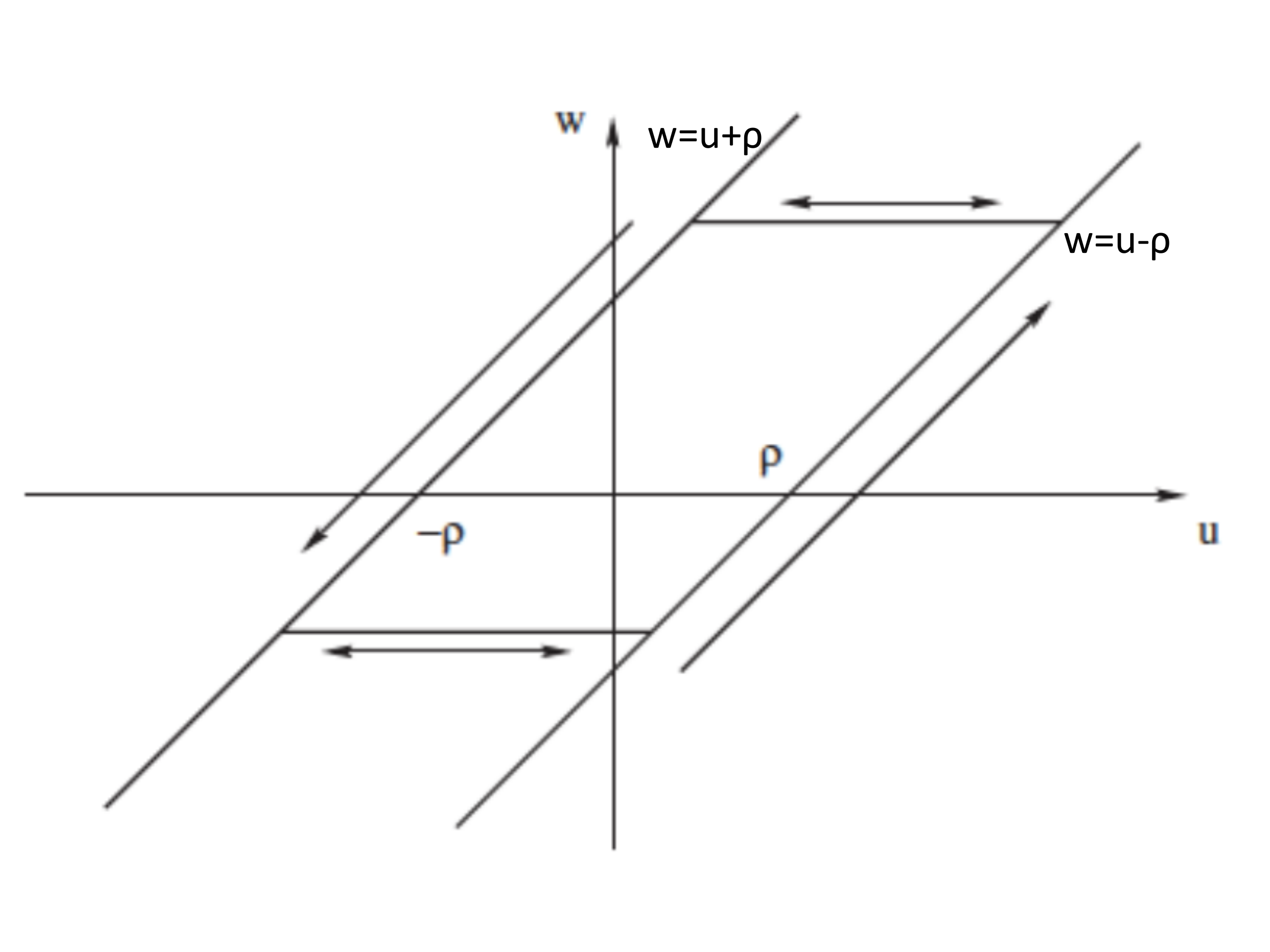}
\caption{\label{Play} Hysteresis play operator }
\end{center}
\end{figure}
In Figure \ref{Play}, $\rho>0$ is the parameter characterizing this operator. We define
\begin{equation}
\label{omega_rho}
\Omega_{\rho}:=\left\{(u, w) \in \mathbb{R}^2 \ \vert u-\rho < w < u+\rho\right\}.
\end{equation}
Given a scalar input $u$ (a continuous function of time), the behavior of the output (a continuous function of time, too) of the scalar play operator $w(\cdot):=\mathcal{P}[u](\cdot)$, with its typical hysteresis loops, can be described using  the phase-portrait in Figure \ref{Play}, representing the trajectories $t\mapsto(u(t),w(t))\in\mathbb{R}^2$. In particular, supposing $u$ piecewise monotone, we have the following. If $(u(t), w(t) \in\Omega_\rho$ for all $t\in I$, with $I$ interval, then $w$ is constant in $I$ (the pair $(u,w)$ moves horizontally, in any of the two possible directions); if $w(t) = u(t) -\rho$ (i.e. $(u,w)$ belongs to the right-boundary of the strip $\Omega_\rho$), $u$ is non increasing in $[t, t + \tau] $ and $w(t)\le u(t+\tau)+\rho$ (i.e. in $[t,t+\tau]$ the pair $(u(\cdot),w(\cdot))$ does not go out the closed strip $\overline\Omega_\rho$) then $w$ stays constant in $[t,t + \tau]$ (again, $(u,w)$ only moves horizontally); if $w(t) = u(t)-\rho$ and $u$ is nondecreasing in $[t,t + \tau]$ then $w = u -\rho$ in $[t,t + \tau]$  (i.e. $(u,w)$ moves along the right-boundary of $\Omega_\rho$ in the upward versus only); a similar argumentation holds if $w(t) = u(t) + \rho$. Moreover we have to prescribe also an initial value for the output: $w(0) = w_0$, with the condition $(u(0),w_0)\in\overline\Omega_\rho$. Finally, we point out the memory feature of the play operator: for a given value of the input, $u(t)$, there is in principle a whole interval of possible values for the output: $[u(t)-\rho,u(t)+\rho]$, and the actual value depends on the past history of the input.

The previous discussion has assumed that the input $u(\cdot)$ is piecewise monotone. However, due to its good continuity/convergence properties, the play operator can be also defined for any continuous input, using an approximation of the input by a sequence of piecewise monotone functions (see Krasnosel'ski\v{\i}-Pokrovski\v{\i} \cite{KR} and Visintin \cite{Visintin94}). In particular, the phase-portrait in Figure \ref{Play} is still preserved. Given a time horizon $T>0$, a possible characterization of the output for an absolutely continuous input $u\in W^{1,1}(0,T)$, with initial output $w_0$ such that $|w_0-u(0)|\leq\rho$, is as the unique absolutely continuous function $w$ such that

\[
\left\{
\begin{array}{ll}
\displaystyle
|u(t)-w(t)|\le\rho&\forall\ t\\
\displaystyle
\frac{dw}{dt}(t)(u(t)-w(t)-v)\ge0&\forall v\ \mbox{such that } |v|\le\rho,\ \mbox{a.e. } t\\
\displaystyle
w(0)=w_0
\end{array}
\right.
\]

Hence, the play operator $\cal P$ is finally defined on the space of the continuous functions, more precisely, in the set

$$
\begin{array}{ll}
\displaystyle
{\mathcal D}:=\left\{(u,w_0)\in C^0([0,T])\times{\mathbb R}|(u(0),w_0)\in\overline\Omega_\rho\right\}\\
\displaystyle
\mathcal{P}:{\mathcal D}\to C^0([0,T]),\ (u,w_0)\mapsto w:={\cal P}[u,w_0]
\end{array}
$$

\noindent
where $ C^0([0,T])$ is the set of all continuous functions defined in $[0,T]$, $T>0$.

The play operator is used in literature to model several hysteresis phenomena (such as, for example, the mechanical play in a junction (also called backlash) due to some damage) and moreover it has many interesting and useful properties (see for example Visintin \cite{Visintin94}): For every $(u,w_0),(v,w_0)\in{\mathcal D}$:
\begin{itemize}
\item[a)] Causality: $ \footnotesize u|_{[0,t]}=v|_{[0,t]} \Rightarrow \mathcal{P}[u,w_0](t)=\mathcal{P}[v,w_0](t)$.
\item[b)] Rate independence: $\mathcal{P}[u\circ \phi, w_0]=\mathcal{P}[u,w_0]\circ\phi$, where $\phi$ is any time re-parametrization, continuous and non decreasing.
\end{itemize}
The two properties above are shared by almost all hysteresis operators. Other useful properties of the play are:

\begin{itemize}
\item[c)] Lipschitz continuity: $\exists L>0$ such that  $\|\mathcal{P}[u,w^1_0]-\mathcal{P}[v,w^2_0]\|_{C^0([0,T])}\leq L(\|u-v\|_{C^0([0,T])}+|w^1_0-w^2_0|),\ \forall\ (u,w_0^1),(v,w_0^2)\in{\mathcal D},\ t\in[0,T]$,
\item[d)] Semigroup property: $\mathcal{P}[u,w_0](t)=\mathcal{P}[u|_{[\tau,t]},\mathcal{P}[u,w_0](\tau)](t-\tau)$ $\forall 0\leq\tau\leq t,\ \forall\ (u,w_0)\in{\mathcal D},\ t\in[0,T]$,
\end{itemize}
where $\|\cdot\|_{C^0([0,T])}$ denotes the uniform norm in $C^0([0,T])$.\\
These facts make the play operator an easy and good model for our purposes, both from a mathematical and applicative point of view. 

\begin{remark}\label{rmrk:discontinuous} In the sequel we will always consider the play operator as applied to continuous inputs. And this is how it is introduced in 
Krasnosel'ski\v{\i}-Pokrovski\v{\i} \cite{KR} and in Visintin \cite{Visintin94}.
Possible extensions to the case of discontinuous inputs have been analyzed in Brokate-Sprekels \cite{brospr}, Krej\v{c}\'i-Lauren\c{c}ot \cite{Kreici}, Recupero  \cite{VR1,VR2}. 
Note that, when considering discontinuous inputs, in particular a special kind of jump/continuous functions (the so-called regulated functions), one has to decide how to fill the gap in the jumps, in order to recover a (at least approximating) continuous input.  But a simple $L^1$-convergence of the continuous approximation is not sufficient for maintaining the memory feature of the operator, because $L^1$-convergence does not detect the peaks of the functions, which are instead very important in the hysteresis, rate/independent, memory effects. Indeed, in that quoted literature, the required convergence takes also account of some kinds of convergence as functions with bounded variations (BV). In general, in control problems, a BV-convergence is too strong, especially when speaking about convergence of controls. 
However, in the case of hysteresis in the controls, Section 2, the two approaches (continuous and discontinuous inputs) may somehow overlap, leading to the same result. 
In the case of hysteresis in the space, Section 3.2 and 3.3, we are instead in some sense forced to use the continuous inputs approach because the discontinuous one would require the use of controls which are atomic measures instead of measurable functions, as we require in our model. Also note that we have mainly in mind mechanical applications of our results for which atomic measure controls are probably not suitable, and moreover their use is certainly forbidden (even if approximated) in the possible case when there is some restrictions on the boundedness of the admissible controls. See also Remarks \ref{rmrk:previous}--\ref{rmrk:bound}.\end{remark}

The second operator on which we focus is the delayed relay, which models a switching input-output relationship between a time continuous scalar input $z$ and a discrete time dependent output $w\in\{1,-1\}$. 
\begin{figure}[H]
\begin{center}
\includegraphics[scale=0.35]{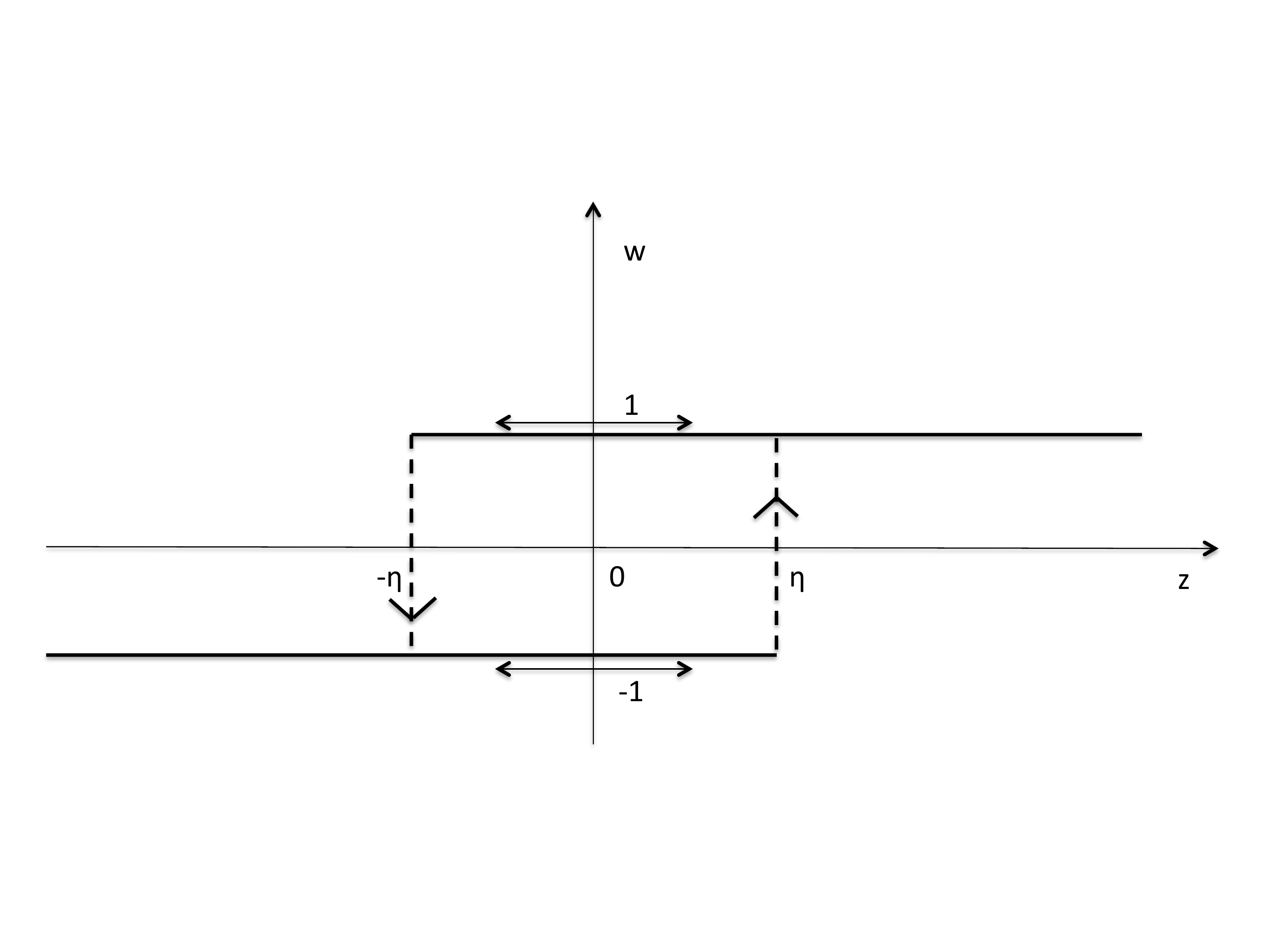}
\caption{\label{Thermostat} The delayed relay operator with thresholds $(-\eta,\eta)$}
\end{center}
\end{figure}
Also in this case, we explain the behavior of the relation $z\mapsto w$ by using Figure \ref{Thermostat} where the phase-portrait of such a delayed switching rule is reported, i.e. the switching trajectories of the pair $(z,w)$ are represented. In particular Figure \ref{Thermostat} corresponds to a delayed swithcing rule with thresholds $(-\eta,\eta)$, $\eta>0$. For example, suppose that at certain time $t$, $w(t)=1$. This means that we certainly have $z(t)\ge-\eta$. The output $w$ will remain equal to $1$, until $z$ will remain larger than or equal to $-\eta$ (i.e. until the pair $(z,w)$ will belong to the closed line $[-\eta,+\infty[\times\{1\}$). If, at a certain time, $z$ becomes strictly lower than $-\eta$, then $w$ switches to $-1$ and it will remain equal to $-1$ until $z$ will possibly become strictly larger than $\eta$.
This is a hysteretic behavior with rate-independent memory: when $-\eta\le z\le\eta$ then the value of $w$ depends on the past history of $z$. Hence we also need an initial value of the output, $w_0\in\{-1,1\}$ such that $(z(0),w_0)\in\left(]-\infty,-\eta]\times\{-1,1\}\right)\cup\left([-\eta,+\infty[\times\{-1,1\}\right)$.
See Visintin \cite{Visintin94} for more details and a possible analytical description of such a behavior. Now denote by $w(\cdot)=h_{\eta}[z](\cdot)$ the delayed relay thermostat, and consider the scalar ODE
$$
\begin{cases}
\dot{z}=g(z,w)\\
w=h_{\eta}[z]\\
z(0)=z_0\quad w(0)=w_0
\end{cases}
$$
with $g$ a suitably regular function. A solution is an absolute continuous function $z(\cdot)$ which solves the ODE in any interval where $w$ is constant and the switching in $w$ occurs when, keeping moving with the previous mode, the solution would be forced to cross the corresponding threshold (see Bagagiolo \cite{Bagagiolo01}). Note that, due to the delayed thresholds, $-\eta<0<\eta$, the solution $z$, in any compact time-interval, can pass to one thresholds to the other just a fixed number of times.

Finally we recall that the play operator can be seen as the superposition of infinitely many delayed relay (see Mayergoyz \cite{Mayer} and Visintin \cite{Visintin94}). We are going to better explain, and use, such a kind of approximation in Subsection \ref{S_P_Hysteresis}.

\section{Hysteresis in the controls}
\label{sec:hysteresis_in_control}
Let us consider a control affine driftless system of the type 
\begin{equation}
\begin{cases}
\label{system_non_hyst}
\dot{\vecz}=\sum_{i=1}^{m}\vecg_i(\vecz)u_i\\
\vecz(0)=\vecz_0
\end{cases}
\end{equation}
where $\vecz(\cdot),\vecz_0\in\mathbb{R}^n$ and the vector fields $\vecg_i:\mathbb{R}^n\to\mathbb{R}^n$ are of class $C^{\infty}$.
\begin{definition}
\begin{itemize}
\item [i)]
We denote by $\mathcal{U}$ (\textbf{admissible controls}) the set of the functions $\vecu=(u_i)_{i=1}^m:[0,+\infty[\to\mathbb{R}^m$ whose components $u_i$ are piecewise constant functions. 
\item [ii)]
The system \eqref{system_non_hyst} is said to be \textbf{controllable } if for any two points $A=\vecz_A$ and $B=\vecz_B$ in $\mathbb{R}^n$ there exists an admissible control $\vecu\in\mathcal{U}$, defined on some time interval $[0, T]$, such that the trajectory of the system \eqref{system_non_hyst}, with initial condition $A$, reaches the point $B$ in time $T$. If for all $A,B\in\mathbb{R}^n$ we can choose $T$ independently from $A$ and $B$, then the system is said to be \textbf{controllable in time $T$.}
\end{itemize}
\end{definition}

Note that, in general, we cannot pretend to control the system using only continuous controls. However, under some suitable hypotheses, the class $\mathcal U$ of piece-wise constant controls, is sufficient. The next assumption goes in that direction.

\begin{hypo}\label{hypo_controllability}
System \eqref{system_non_hyst} satisfies the so-called Chow hypothesis (see Coron \cite{Coron56}), more precisely
the Lie algebra generated by the vector fields $\vecg_i$ is fully generated, i.e. 
$$
\dim\Bigl(\mathscr{L}\text{ie}\{\vecg_i,\,i=1\cdots m\}\Bigr)=n
$$
where $\mathscr{L}\text{ie}\{\vecg_i,\,i=1\cdots m\}$ is the space of the linear combinations $X$ of iterated Lie brackets of the vector fields: $X=\sum_{\ell=1}^q\lambda_{\ell}Y_{\ell}$ with $Y_{\ell}=[\vecg_{k},[\cdots,[\ldots,[\vecg_{i},\vecg_{j}]\ldots]]]$ for $i,j,k=1\ldots m$, and $[\vecg_i,\vecg_j]=\nabla \vecg_j\cdot \vecg_i-\nabla \vecg_i\cdot \vecg_j$
\end{hypo}
In particular, Assumption \ref{hypo_controllability} guarantees that the system is controllable in time $T$ for all $T$.


Let us now consider the following nonlinear system:
\begin{equation}
\begin{cases}
\label{system_hyst_control}
\dot{\vecz}=\sum_{i=1}^{m}\vecg_i(\vecz)\mathcal{P}[v_i,w_0^i]\\
\vecz(0)=\vecz_0
\end{cases}
\end{equation}
where $\mathcal{P}[v_i,w_0^i]$ is the play operator applied to the input $v_i$, with $(v_i(0),w_0^i)\in\overline\Omega_\rho$, and the inputs $v_i$ are at disposal of the controller. System \eqref{system_hyst_control} is obtained from system \eqref{system_non_hyst} by replacing any control $u_i$ iwith $\mathcal{P}[v_i,w_0^i]$. If $u_i=v_i$, then it means that system \eqref{system_non_hyst} is not directly experiencing the actuation of the control $u_i$, but instead a sort of perturbation of it. This can be due, for example, to some kinds of damage in the mechanical realization of the control or to some kinds of general hysteresis effect: think of the case where $\vecu=(u_i)_i$ is the magnetic field but the system reacts to the magnetization of a ferromagnetic actuator, represented here by the output of the play operator. However, note that in this case both the actual controls $\mathcal{P}[v_i,w_0^i]$ and their inputs $v_i$ are not in $\mathcal{U}$ (piecewise constant functions) but, due to the construction of the play operator, must belong to the space of continuous functions $C^0$.

Our goal is to investigate the controllability properties of system \eqref{system_hyst_control}. We are going to use the following result (here and in the sequel for an interval $I$, $\chi_I$ is its characteristic function: $\chi_I(t)=1$ if $t\in I$ $\chi_I(t)=0$ otherwise.)
\begin{lemma}
\label{approx_controls}
For every piecewise constant function $\bar{\vecu}:[0,T]\to\mathbb{R}^m$ and for every initial state $\vecw_0=(w_0^{i})_{i=1}^m$, there exists a sequence of continuous functions $\left(\vecv^k=(v_i^{k})_{i=1}^m\right)_{k\in\mathbb{N}}$, such that $\vecu^k:=\mathcal{P}[\vecv^k,\vecw_0]$ converges to  $\bar{\vecu}$ in $L^1(0,T)$ as $k\to+\infty$.
Here $\mathcal{P}[\vecv^k,\vecw_0]$ stays for the vector  $\left(\mathcal{P}[v^k_i,w_0^i]\right)_{i=1}^m$.
\end{lemma}
\Proof
First of all note that the components $\bar{u}_i$ of $\bar{\vecu}$ are not continuous (piecewise constant), therefore, we cannot exactly generate them as outputs of the play operator, since the output of the play is a continuous function. Hence we first approximate each $\bar{u}_i$ with a sequence of continuous piecewise linear functions. 

Of course, we can argue for any single scalar component $\overline u_i$ and hence we drop the notation of the index $i$ for simplicity. Let $\bar{u}(t)=\sum_{j=1}^n\alpha_j\chi_{[t_{j-1},t_j]}(t)$ be a piece-wise constant function (with $0=t_0<t_1<\cdots t_n=T$, and $\alpha_j\in\mathbb{R}$). For every integer $k>0$ sufficiently large,  we consider the piece-wise linear function $u^k$ defined by

\begin{equation}
\begin{aligned}
\label{uk}
\small
u^k(t)=&\chi_{[0,\frac{1}{k}]}(t)\Bigl(kt(\alpha_1-w_0)+w_0\Bigr)+\alpha_1\chi_{[\frac{1}{k},t_1-\frac{1}{k}]}(t)+\\
&+\sum_{j=1}^{n-1}\Bigl((\alpha_{j+1}-\alpha_j)\bigl(\frac{kt}{2}-\frac{k}{2}(t_j-\frac{1}{k})\bigr)+\alpha_j \Bigr)\chi_{[t_j-\frac{1}{k},t_j+\frac{1}{k}]}(t)+\\
&+\sum_{j=1}^{n-1}\alpha_j\chi_{[t_j+\frac{1}{k},t_{j+1}-\frac{1}{k}]}+\alpha_n\chi_{[t_{n-1}+\frac{1}{k},t_n]}(t)
\end{aligned}
\end{equation}

\noindent
\begin{figure}[H]
\begin{center}
\includegraphics[scale=0.45]{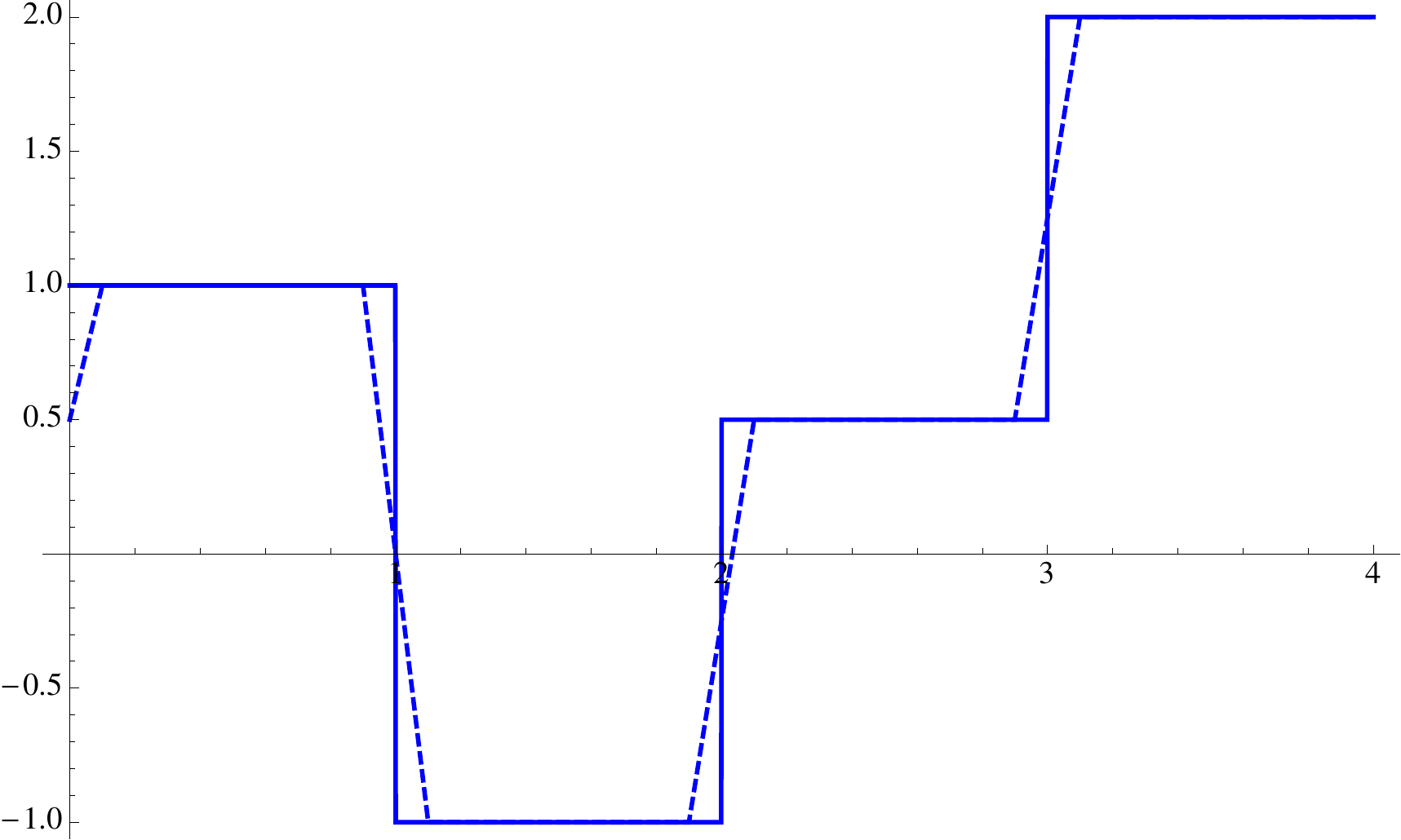}
\caption{\label{control_approx} The solid line is the piecewise constant control $\bar{u}$ the dashed one, one of the approximating $u^k$ (see \eqref{uk}). Here it is $T=4$, $n=4$, $t_1=1,\,t_2=2\,t_3=3,\,t_4=4$ $\alpha_1=1\,\alpha_2=-1,\,\alpha_3=\frac{1}{2},\alpha_4=2$, $k=10$ and $w_0=\frac{1}{2}$}
\end{center}
\end{figure}

\begin{figure}[H]
\begin{center}
\includegraphics[scale=0.3]{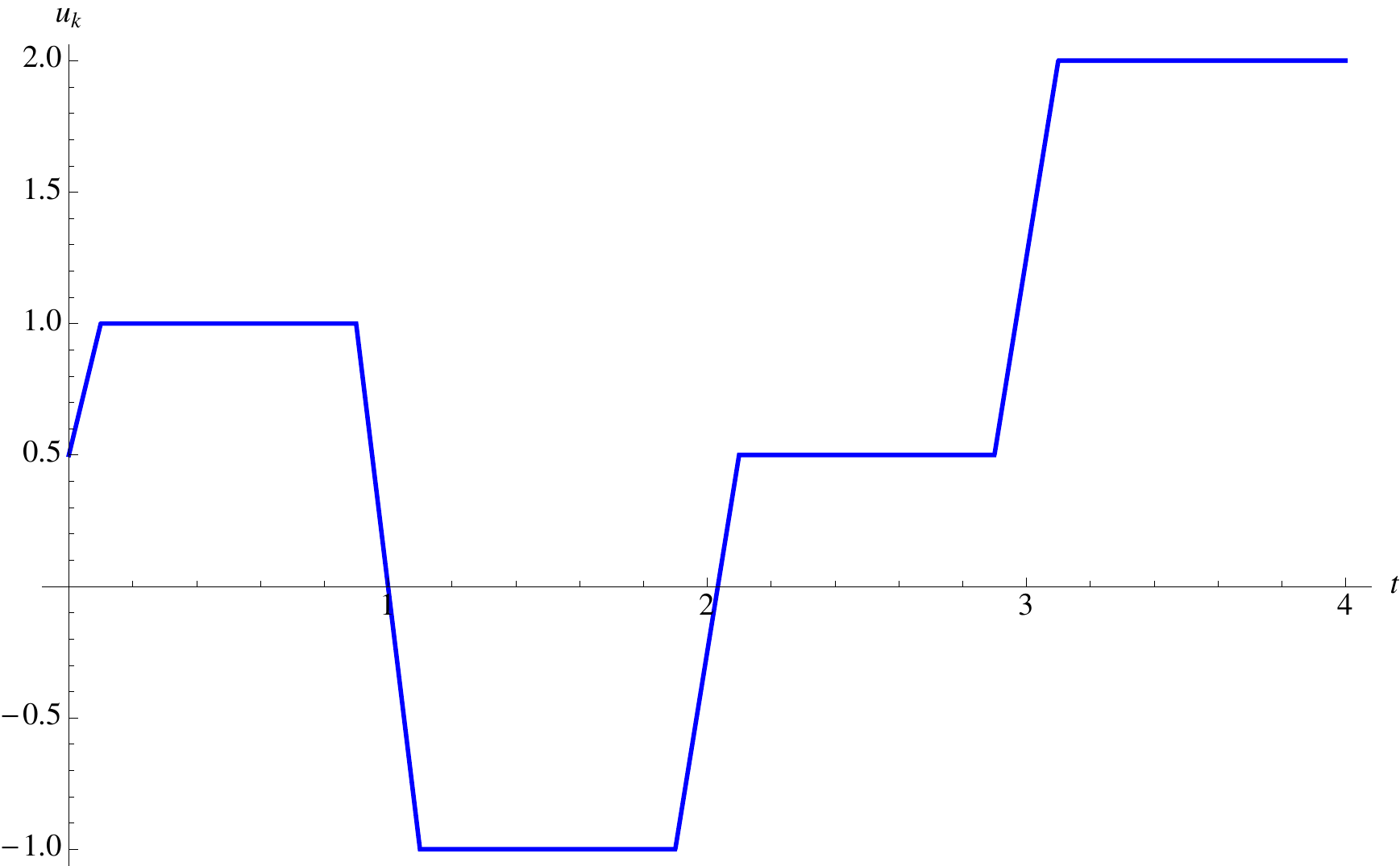}
\includegraphics[scale=0.3]{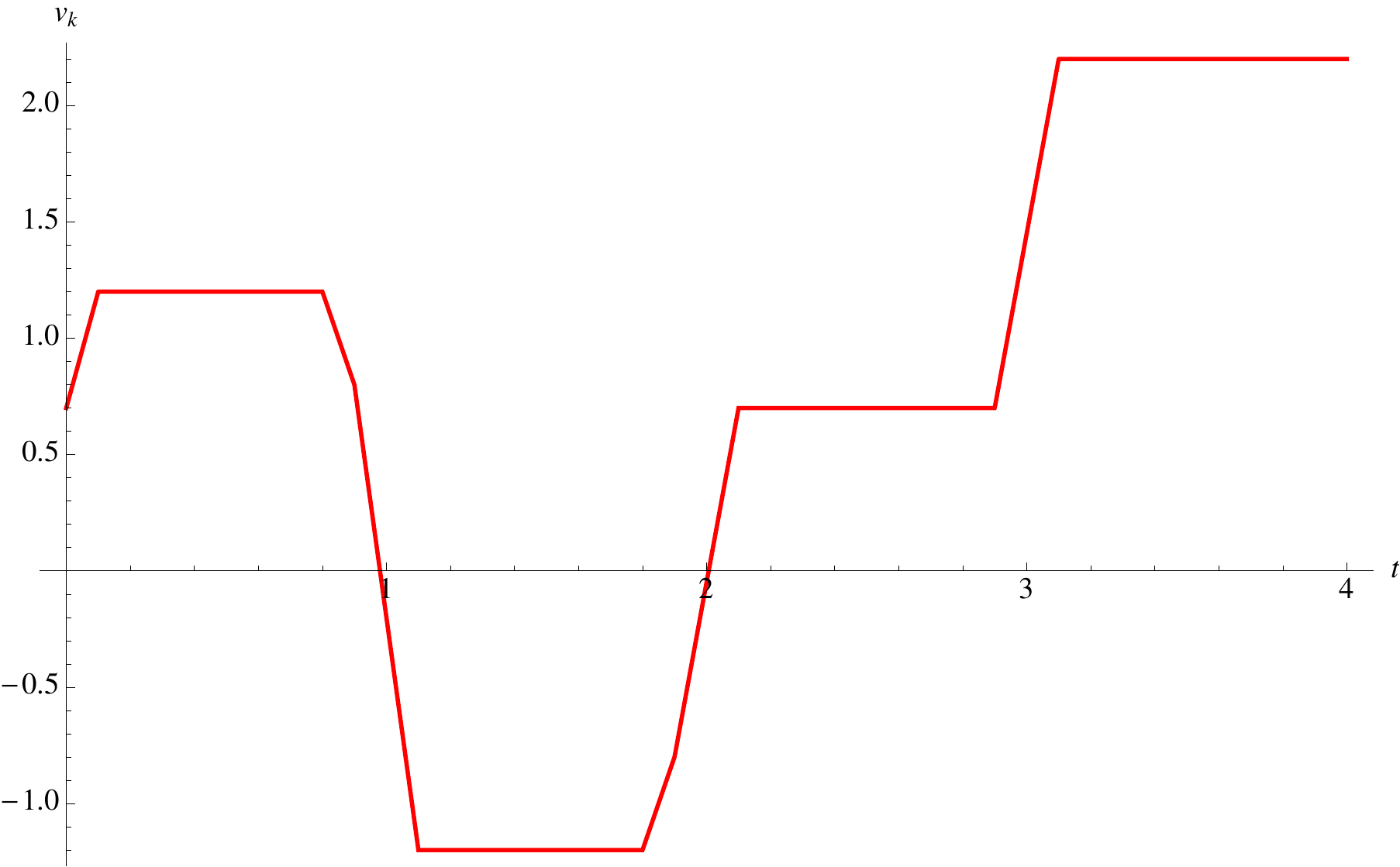}
\caption{\label{play_approx} On the left the function \CorrR{$u^k$ defined in \eqref{uk} with the values as in Figure \ref{control_approx}, on the right the corresponding input of the play operator  $v^k$ defined in \eqref{vk}, for $\rho=0.2$}. }
\end{center}
\end{figure}

\noindent
In particular note that $u^k(0)=w_0$.
By construction, the convergence $u_{k}\to\bar{u}$ in $L^{1}(0,T)$ is immediate. 
For any $k$ we then consider the following piece-wise linear function

\begin{equation}
\begin{aligned}
\label{vk}
&v^k(t)=\chi_{[0,\frac{1}{k}]}(t)\bigl(u^k(t)+sgn(\alpha_1-w_0)\rho\bigr)+(u^k(t)+sgn(\alpha_{2}-\alpha_1)\rho)\chi_{[\frac{1}{k},t_1-\frac{2}{k}]}+\\
&+\sum_{i=1}^{n-2}\Bigl(\bigl(u^k(t)+sgn(\alpha_{i+1}-\alpha_i)\rho\bigr)\chi_{[t_i-\frac{2}{k},t_{i+1}-\frac{2}{k}]}\Bigr)\Theta((\alpha_{i+1}-\alpha_i)(\alpha_{i+2}-\alpha_{i+1}))\\
&+\Bigl(\bigl(sign(\alpha_{i+1}-\alpha_i)2\rho k( t-t_i+\frac{2}{k})+\alpha_i-sign(\alpha_{i+1}-\alpha_i)\bigr)\chi_{[t_i-\frac{2}{k},t_i-\frac{1}{k}]}+\\
&(u^k(t)+sign(\alpha_{i+1}-\alpha_i))\chi_{[t_i-\frac{1}{k},t_{i+1}-\frac{2}{k}]}\Bigr)\Theta(-(\alpha_{i+1}-\alpha_i)(\alpha_{i+2}-\alpha_{i+1}))+\\
&+(u^k(t)+ sgn(\alpha_{n}-\alpha_{n-1})\rho)\chi_{[t_{n-1}-\frac{2}{k},t_n]}
\end{aligned}
\end{equation}
where $\Theta$ is the Heaviside function that is $1$ if its argument is positive and $0$ otherwise, and $sgn$ is the sign function. Some calculations may show that $u^k=\mathcal{P}[v^k,w_0]$. In particular, note that, being $u^k(0)=w_0$, we have $(v^k(0),w_0)\in\overline\Omega_\rho$ and also that $v^k$ is piece-wise monotone. Hence we can easily perform the constructive description of the output of the play as in the previous section. Looking to the example to Figure \ref{play_approx}, $v^k(0)=u^k(0)+\rho$ and, in the interval $[0,1/k]$ $v^k$ and $u^k$ increase together, reaching the values $u(t/k)=\alpha_1, v(t/k)=\alpha_1+\rho$; in the interval $[1/k,t_1-2/k]$ they both remain constant; in the interval $[t_1-2/k,t_1-1/k]$ $v^k$ rapidly decreases to the value $\alpha_1-\rho$, $u^k$ remaining constant; in the interval $[t_1-1/k,t_1+1/k]$ they both decrease together, reaching the values $u^k(t_1+1/k)=\alpha_2, v^k(t_1+1/k)=\alpha_2-\rho$. We proceed in this way, concluding the proof. \EndProof
\begin{remark} \label{surjectivity}
Observe that Lemma \ref{approx_controls} implies that on the class of functions of the type \eqref{uk}, the play operator is surjective, since there always exists an input like \eqref{vk} whose output through the play operator is exactly \eqref{uk}.
\end{remark}

\begin{remark} Note that $u^k$ and $v^k$ have sometimes to change their values by some fixed amplitudes ($2\rho$, or $\alpha_{i+1}-\alpha_i$) in a time interval of length $1/k$. Hence their derivatives diverge when $k\to+\infty$. However, this is not a problem for our construction because, the convergence of $u^k$ to $\overline u$ is only required in $L^1$, and the derivatives of $v^k$ do not play any role, due to the rate-independence of the play operator: only the sequence of values reached by $v^k$ in its history has a role.
\end{remark}

\begin{theorem}

\label{approx_control}
Let us assume the controllability of the non hysteretic system \eqref{system_non_hyst}, Assumption \ref{hypo_controllability}. Given two points $A,B\in\mathbb{R}^n$ and $T>0$, let us consider the piece-wise constant control $\overline\vecu\in{\mathcal U}$ which steers system \eqref{system_non_hyst} from $A$ to $B$ in time $T$, and let $\vecz$ be the corresponding trajectory. We consider the corresponding controls $\vecv^k$ 
defined as in Lemma \ref{approx_controls}. 
Then for any given initial data for the play operator, $\vecw_0$, the sequence of the trajectories $\vecz_k$ of the systems
\begin{equation}
\begin{cases}
\label{system_hyst}
\dot{\vecz}_k=\sum_{i=1}^{m}\vecg_i(\vecz_k)\mathcal{P}[v^k_i,w_0^i]\\
\vecz_k(0)=A
\end{cases}
\end{equation}
uniformly converges on $[0,T]$, as $k\to+\infty$, to the trajectory $\vecz$ of the non hysteretic system \eqref{system_non_hyst} with controls $\bar{\vecu}$.
\end{theorem}
\Proof
Let  $\vecz_k$ be the unique solution of \eqref{system_hyst}. 
Denoting by $\vecu^k$ the output of the play operator in \eqref{system_hyst}, by Lemma \ref{approx_controls} $\vecu^k\to\bar{\vecu}$ in $L^1(0,T)$ componentwise, and $\vecu^k$  are equibounded in $L^{\infty}(0,T)$ by construction. Moreover $\vecg_i\in C^{\infty}$ for all $i$ and thus the solutions $\vecz_k$ are both locally equi-bounded and locally equi-Lipschitz continuous. By hypothesis, there exists a non hysteretic trajectory $\vecz$ of \eqref{system_non_hyst} with controls $\bar{\vecu}$ defined on the whole interval $[0,T]$ (the controlled one, from $A$ to $B$). This implies, by standard estimate arguments on the trajectories, that the solutions $\vecz_k$ of \eqref{system_hyst} exist on the whole $[0,T]$ too, and that they do not exit from a common compact set. 
Indeed, we take a ball $\mathcal{B}$ that contains in its interior the trajectory $\vecz$, and, for every $k$, take $t_k>0$ the possible first instant such that $\vecz_k(t_k)\in\partial B$. Obviously, every $\vecz_k$ is defined at least in $[0,t_k]$. If, by contradiction, there exists a convergent subsequence $t_k\to\overline t\le T$, then, being all contained in $\mathcal B$, where the vectors $\vecg_i$ are equi-bounded (and so are the controls $\vecu^k$) and equi-Lipschitz, reasoning as in the following estimates, we would get the contradiction $\partial \mathcal B\in \vecz_k(t_k)\to\vecz(\overline t)\in{\rm int}\mathcal B$.

Hence in the sequel we will treat the vector fields $\vecg_i$ as uniformly bounded and Lipschitz along the trajectories.
We have the following estimates for all $t\in[0,T]$ (writing $u^k_i=\mathcal{P}[v^k_i,w_0^i]$, and recalling that $\vecz_k(0)=\vecz(0)=A$)
\begin{equation*}
\small
\begin{aligned}
&\|\vecz_k(t)-\vecz(t)\|\leq\\&\sum_{i=1}^m\int_0^t \Bigl(\|\vecg_i(\vecz_k(\cdot))\|_{\infty}|u^k_i(s)-\bar{u}_i(s)|+\|\vecg_i(\vecz_k(s))-\vecg_i(\vecz(s))\|\|\bar{\vecu}\|_{\infty}\Bigr)\,ds
\end{aligned}
\end{equation*}

\noindent
Denoting respectively by $M'$ and $L$ the uniform bound and  Lipschitz constant of the functions $\vecg_i$ and taking $M=max\{M',\|\bar{\vecu}\|_\infty\}$, we get
$$
\|\vecz_k(t)-\vecz(t)\|\leq C_k+m M L \int_0^t\|\vecz_k(s)-\vecz(s)\|ds
$$
where $C_k=m M'\int_0^T |u^k_i(s)-\bar{u}_i(s)|\,ds\rightarrow 0$ for the convergence of $u^k_i$ to $\bar{u}_i$ in $L^1$. The last inequality, for the arbitrariness of $t\in[0,T]$, by  the Gronwall lemma implies
$$
\|\vecz_k-\vecz\|_{C^0([0,T])}\leq C_k e^{mMLT}
$$
which converges to zero as $k$ tends to infinity.
\EndProof

\begin{remark}
Theorem \ref{approx_control} is an approximate controllability result for the hysteretic system \eqref{system_hyst}, i.e. given any $T>0$ any initial and final configurations $A$ and $B$, and any initial output $\vecw_0$, we are always able to find a sequence of piecewise linear controls $\vecv$ which allows us to arrive as close as we want to $B$, but in general, it is not guaranteed that we can reach it.
\end{remark}

\begin{remark}
\label{Robustness}
Note that the construction of the approximating controls $\vecv^k$ obviously depends on the amplitude $\rho>0$ of the considered play operator, which here we denote by ${\mathcal P}_\rho$ (see Figure \ref{Play}). On the other hand, since the construction of the sequence $\left(u^k_i\right)_k$ is independent on the parameter $\rho$, we easily have a sort of robustness of our approximating procedure, in the sense that  $\vecv^k\to\vecu^k$  and $\mathcal{P}_\rho[\vecu^k,\vecw_0]\to\vecu^k$ as $\rho\to 0$, uniformly in time.
\end{remark}
\section{Hysteresis in the state}
\label{sec:hysteresis_in_state}
In this section we start from the same smooth controllable system \eqref{system_non_hyst}, and we want to analyze what happens to its controllability properties when the play hysteresis operator is applied to the state variables (see (\ref{eq:intro})-right). We will focus on a system \eqref{system_non_hyst} with a particular ``triangular" structure, being a more general situation far to be clarified. 
More precisely, we consider the following system of the type of \eqref{system_non_hyst} in $\mathbb{R}^3$

\begin{equation}
\label{general_heis}
\begin{pmatrix}
\dot x\\
\dot y\\
\dot z
\end{pmatrix}=\begin{pmatrix}1\\0\\0\end{pmatrix}u_1+\begin{pmatrix}0\\1\\f(x)\end{pmatrix}u_2
\end{equation}
i.e. $\vecg_1=\begin{pmatrix}1\\0\\0\end{pmatrix}$, $\vecg_2=\begin{pmatrix}0\\1\\f(x)\end{pmatrix}$,  $f\in C^{\infty}$. 

Moreover we suppose that Assumption \ref{hypo_controllability} is satisfied, i.e. the Lie algebra generated by $\vecg_1$ and $\vecg_2$ is fully generated so that the system is controllable. More precisely, since $[\vecg_1,\vecg_2]=\begin{pmatrix}0\\0\\f^{'}(x)\end{pmatrix}$, if $f^{'}\neq0$ then it is
$$
\dim\Bigl(\displaystyle{span}\{\vecg_1,\vecg_2,[\vecg_1,\vecg_2]\}\Bigr)=3
$$

We first give a motivating example of mechanical system with that kind of structure.

\subsection{Example}

System (\ref{general_heis}) is a generalization of the Heisenberg flywheel, (see Montgomery \cite{Montgomery}). A point mass $m$ is constrained to slide along a massless rod connected to a flywheel with moment of inertia $I$, and it is able to rotate about it. Moreover the flywheel is attached to a table by a joint on which it spins freely. This joint is frictionless thus it does not exert any torque on the system. 

\begin{figure}[H]
\begin{center}
\includegraphics[scale=0.40]{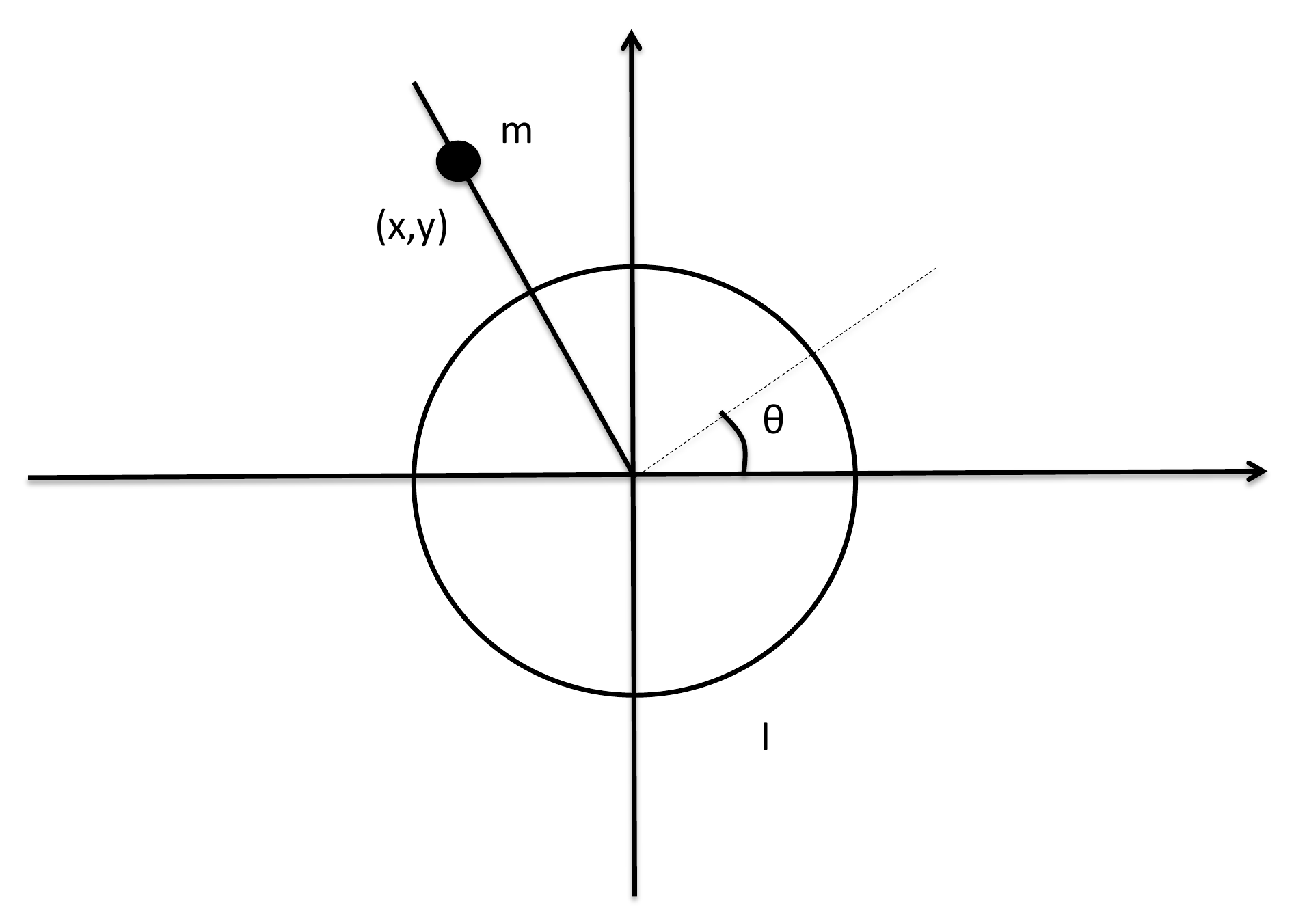}
\caption{\label{Heisenberg_example} The system of the Heisenberg flywheel }
\end{center}
\end{figure}

We denote by $\theta$ the angle of the flywheel relative to the table and by $(x,y)$ the mass coordinates, measured with respect to an external frame. We can exert a torque on the rod to rotate it relatively to the wheel and we are able to slide the mass back and forth on the rod. Therefore we have two controls, the torque $\tau$ and the sliding speed and three states $(x,y,\theta)$. Applying a linear transformation of the controls, the control laws of the system become

\begin{equation}
\label{flywheel}
\begin{pmatrix}
\dot x\\
\dot y\\
\dot \theta
\end{pmatrix}=\begin{pmatrix}1\\0\\\alpha y \end{pmatrix}u_1+\begin{pmatrix}0\\1\\-\alpha x\end{pmatrix}u_2
\end{equation}
 with $\alpha=-\frac{I}{m}$.\\
The vector fields of this control system are exactly the Heisenberg group vector fields. The system \eqref{flywheel} can be rewritten, after the change of coordinates $(x,y,\theta)\to(x,y,z)=(x,y,-\theta/2\alpha+xy/2)$, as the following "triangular" one
\begin{equation}
\label{heis}
\begin{pmatrix}
\dot x\\
\dot y\\
\dot z
\end{pmatrix}=\begin{pmatrix}1\\0\\0\end{pmatrix}u_1+\begin{pmatrix}0\\1\\x\end{pmatrix}u_2
\end{equation}
This system has a lot of good properties, in particular it is controllable. Indeed the Lie algebra generated by the two dynamic vector fields, $\begin{pmatrix}1\\0\\0\end{pmatrix}$ and $\begin{pmatrix}0\\1\\x\end{pmatrix}$ is fully generated. Indeed 
\begin{equation}
\label{lie_heis}
\Bigg[\begin{pmatrix}1\\0\\0\end{pmatrix}, \begin{pmatrix}0\\1\\x\end{pmatrix}\Bigg ]=\begin{pmatrix}0\\0\\1\end{pmatrix} \quad\Rightarrow\quad
\dim\Bigl(\displaystyle{span}\{\begin{pmatrix}1\\0\\0\end{pmatrix}, \begin{pmatrix}0\\1\\x\end{pmatrix},\begin{pmatrix}0\\0\\1\end{pmatrix}\}\Bigr)=3
\end{equation}

\subsection{The hysteretic system}
We are interested in studying the controllability properties of the hysteretic version of (\ref{general_heis})
\begin{equation}
\label{general_heis_hyst}
\begin{pmatrix}
\dot x\\
\dot y\\
\dot z
\end{pmatrix}=\begin{pmatrix}1\\0\\0\end{pmatrix}u_1+\begin{pmatrix}0\\1\\f(\mathcal{P}[x,w_0])\end{pmatrix}u_2
\end{equation}
Note that existence of the trajectory is guaranteed by the Lipschitz property of $\mathcal{P}$.
The following lemma will be used.
\begin{lemma}
\label{density}
The play operator has dense image in the space of continuous piecewise linear functions $x(\cdot)$ in $[0,T]$. That is, for any such $x(\cdot)$, denoting $x_0=x(0)$, there exists a sequence of continuous piecewise linear functions $v^j(\cdot)$ such that $\mathcal{P}[v^j,x_0](\cdot)\to x(\cdot)$ in $L^{\infty}(0,T)$.
\end{lemma}
\Proof
Let us consider a general continuous piecewise linear function
\begin{equation}
\label{x}
x(t)=\sum_{i=0}^{n-1}\bigl(\alpha_{i+1}t+\sum_{k=1}^it_k(\alpha_k-\alpha_{k+1})+x_0\bigr)\chi_{[t_i,t_{i+1}]}(t)
\end{equation}
where $(0=t_0<t_1<\cdots<t_n=T)$ is a subdivision of the interval $[0,T]$ and $\alpha_1,...,\alpha_n\in\mathbb{R}$ are the slopes (with the convention: $\sum_{k=1}^0\xi_k=0$, whichever the quantity $\xi_k$ are). Then consider the following sequence of continuous piecewise linear functions
\begin{equation}
\small
\label{vj}
\begin{aligned}
&v^j(t)=\chi_{[0,\frac{1}{j}]}(t)\Bigl[\alpha_1t+sgn(\alpha_1)\rho+x_0\Bigr]+\sum_{i=0}^{n-1}\Bigl\{\chi_{[t_i+\frac{1}{j},t_{i+1}-\frac{1}{j}]}(t)(x(t)+sgn(\alpha_{i+1})\rho)+\\
&+\Theta(\alpha_{i+1}\alpha_{i+2})\chi_{[t_{i+1}-\frac{1}{j},t_{i+1}+\frac{1}{j}]}(t)(x(t)+sgn(\alpha_{i+1})\rho)+\\
&+\Theta(-\alpha_{i+1}\alpha_{i+2})\chi_{[t_{i+1}-\frac{1}{j},t_{i+1}+\frac{1}{j}]}(t)\Bigr[\frac{j}{2}\Bigl(\frac{1}{j}(\alpha_{i+1}+\alpha_{i+2})+2sgn(\alpha_{i+2}-\alpha_{i+1})\rho\Bigr)t+\\
&+\sum_{k=1}^it_k(\alpha_k-\alpha_{k+1})+x_0+sgn(\alpha_{i+1})\rho\Bigr]\Bigr\}
\end{aligned}
\end{equation}
where $\Theta(\cdot)$ is the Heaveside function.


\begin{figure}[H]
\begin{center}
\includegraphics[scale=0.3]{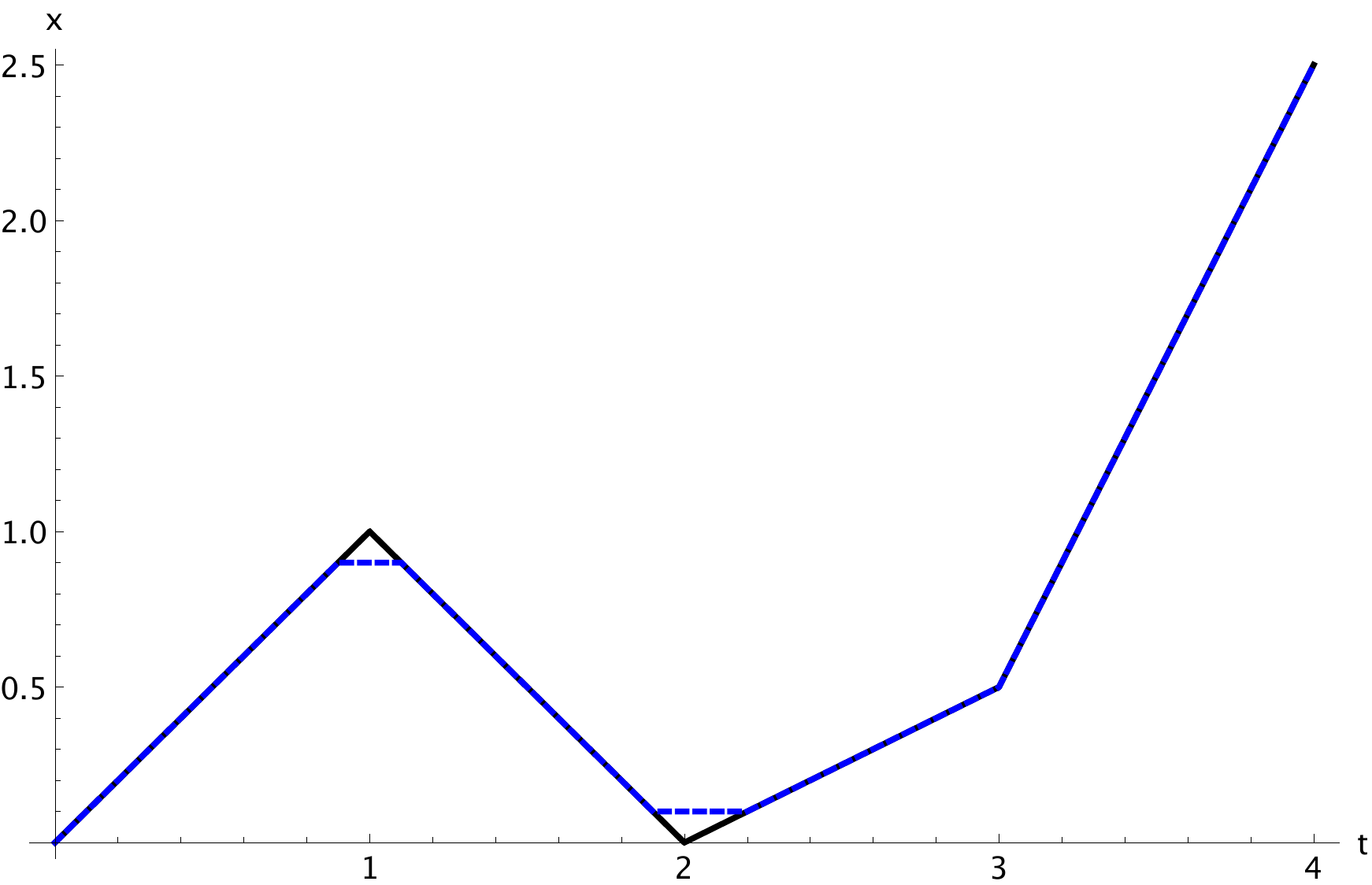}
\includegraphics[scale=0.3]{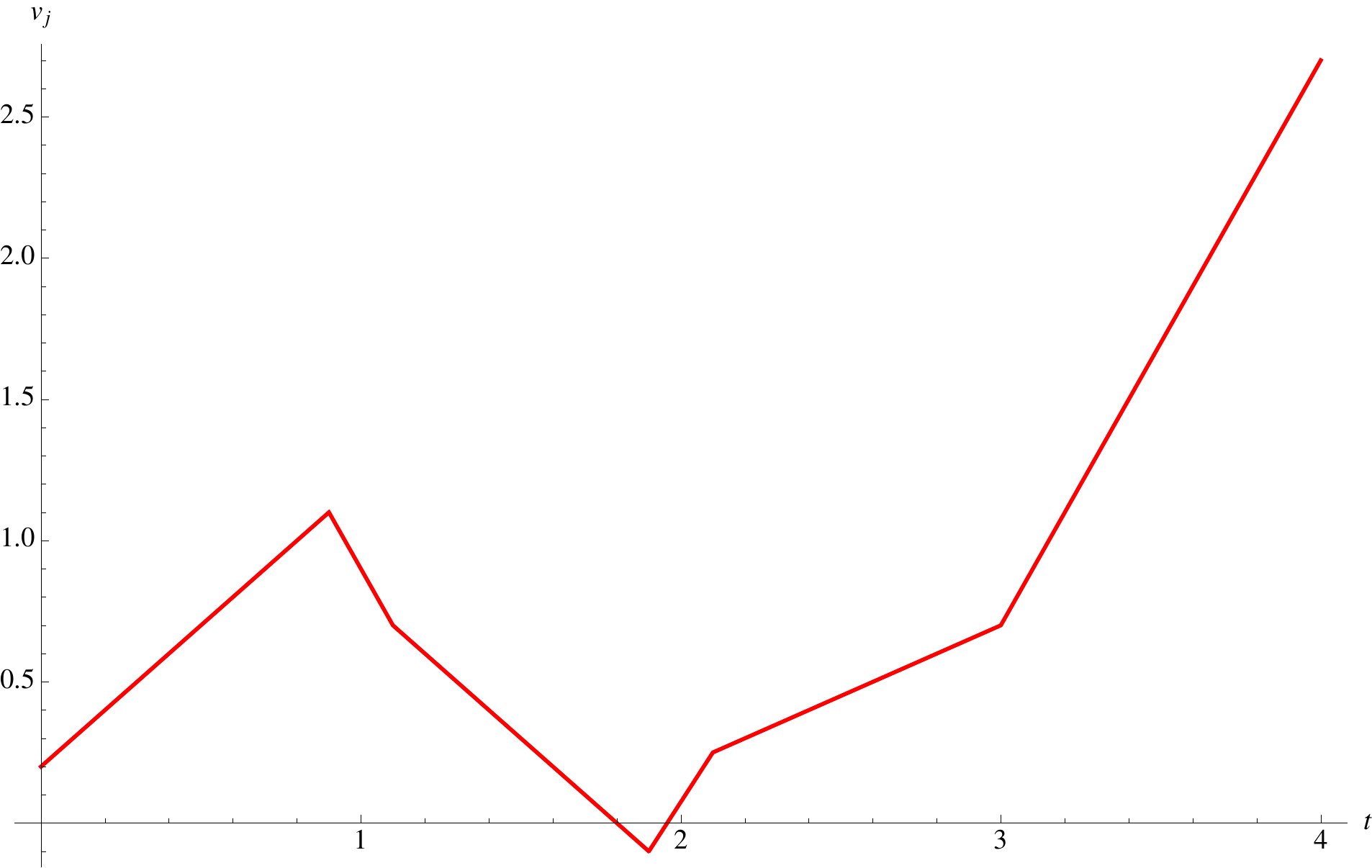}
\caption{\label{traj_approx} On the left, \CorrR{the solid line} represents the  continuous piecewise linear function $x(\cdot)$ defined in \eqref{x} and \CorrR{the dashed line} the approximating output $\mathcal{P}[v_j,x_0]$ with $x_0=0$, $n=4$, $t_1=1,\,t_2=2,\,t_3=3,\,t_4=4$, $\alpha_1=1,\,\alpha_2=-1,\alpha_3=\frac{1}{2},\,\alpha_4=2$ and $j=10$, on the right the input of the play operator with $\rho=0.2$,   $v^j(\cdot)$ defined in \eqref{vj}.}
\end{center}
\end{figure}
\noindent
Note that, applying the play operator to \eqref{vj}, then in the time intervals $[t_i+\frac{1}{j}, t_{i+1}-\frac{1}{j}] $ it is $P[v^j,x_0]=x$. Instead in the time intervals $[t_{i+1}-\frac{1}{j},t_{i+1}+\frac{1}{j}]$, it is $||P[v^j,x_0]-x||_{\infty}= \frac{1}{j}|\alpha_{i+1}| $ which goes to zero when $j\to\infty$. We then have
$$
\mathcal{P}[v^j,x_0](\cdot)\to x(\cdot)\quad\text{for $j\to\infty$ in $L^{\infty}(0,T)$}
$$
\EndProof

Let us suppose that system \eqref{general_heis} (the non hysteretic one) is controllable, for example as said before suppose that the Lie algebra of vector fields $\vecg_1$ and $\vecg_2$ is fully generated. Then for any initial and final conditions $A=(x_A,y_A,z_A)$ and $B=(x_B,y_B,z_B)$, there exists $(\bar{u}_1,\bar{u_2})$ piecewise constant that steers system (\ref{general_heis}) from $A$ to $B$ in time $T>0$.
In particular, due to the structure of the vector fields, $\bar{u}_1$ generates a corresponding continuous piecewise linear trajectory $\bar{x}$ s.t $\bar{x}(0)=x_A,\, \bar{x}(T)=x_B$, moreover  $\bar{u}_2$ is such that, together with $f(\bar{x})$, it generates trajectories $\bar y,\bar z$ such that $\bar{y}(0)=y_A,\, \bar{y}(T)=y_B$ and $\bar{z}(0)=z_A,\, \bar{z}(T)=z_B$.

To make system \eqref{general_heis_hyst} (the hysteretic one) to move from $A$ to $B$, the idea is to look for a control $u_1(\cdot)$ whose integral $x(\cdot)$ is such that $\bar{x}(\cdot)=\mathcal{P}[x,w_0](\cdot)$. In this way using such a control $u_1$ and the same $\bar{u}_2$, as in the non hysteretic case, we can steer $(x_A,y_A,z_A)$ to a point whose second coordinates are the desired ones $(y_B,z_B)$, letting the tuning of the first one to a later time. More precisely we have the following result.
\begin{theorem}
\label{thm3}
Given the controllability of the non hysteretic system (\ref{general_heis}), for any initial and final configurations $A$ and $B$ and for any $w_0$ such that $(x_A,w_0)\in\bar{\Omega}_\rho$, 
there exists a sequence of piecewise constant controls $(u_1^j,u_2^j)$ and a final time $T^*\ge T$ such that the solution $(x^j(\cdot),y^j(\cdot),z^j(\cdot))$ of system \eqref{general_heis_hyst} starting from $A$ is such that
\begin{align*}
&x^j(T^*)= x_B&&y^j(T^*)= y_B&&z^j(T^*)\to z_B\quad\text{as $j\to\infty$}
\end{align*}
\end{theorem}

\begin{figure}[H]
\begin{center}
\includegraphics[scale=0.40]{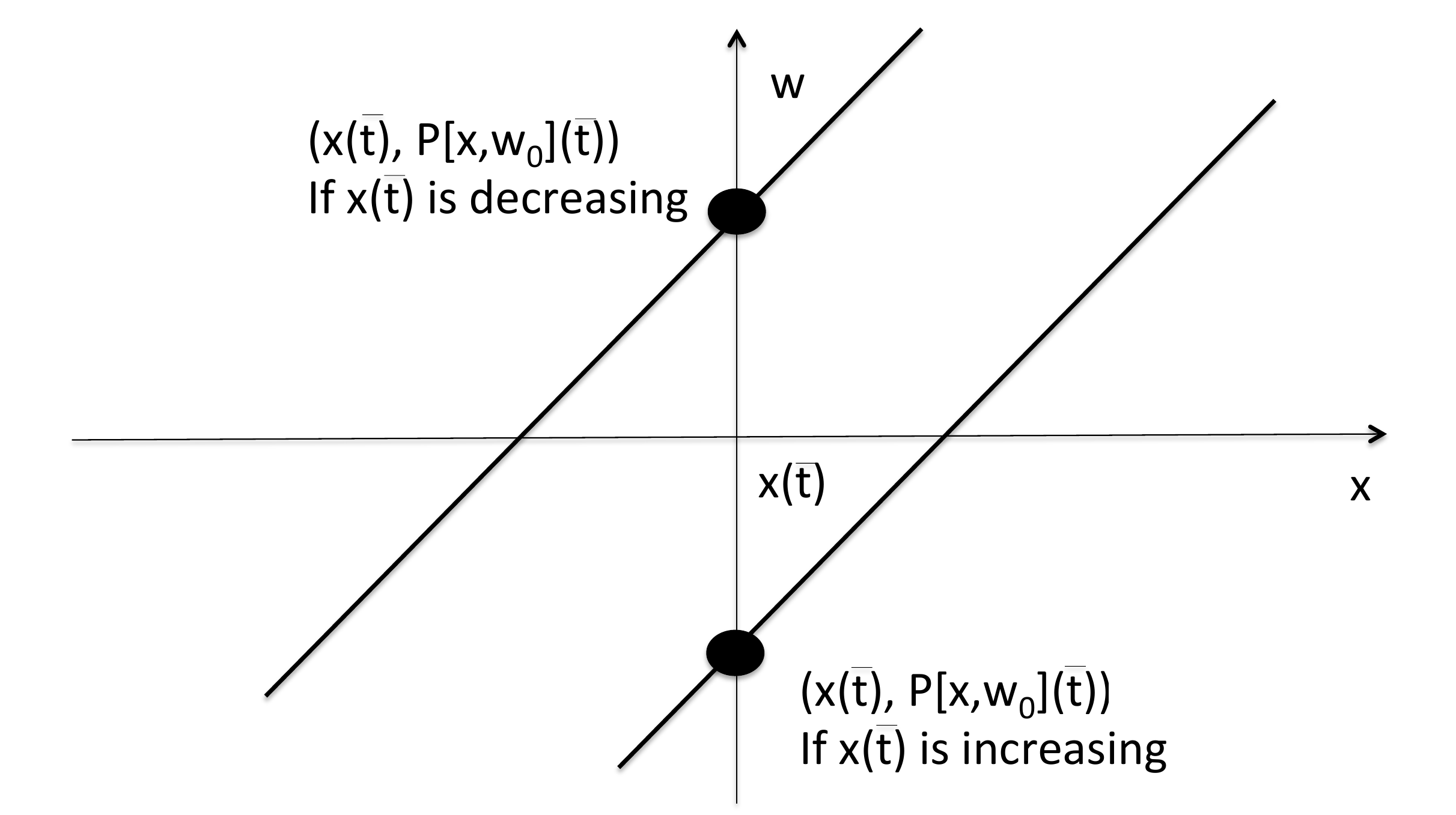}
\caption{\label{initial_data}  From this picture it is evident that the initial couple $(x_0,w_0)$ should stay on the line. }
\end{center}
\end{figure}

\Proof
Let $\bar{\vecu}=(\overline u_1,\overline u_2)$ be the piecewise constant control which steers \eqref{general_heis} from $A$ to $B$, and let $\bar{x}$  be the corresponding continuous piecewise linear trajectory such that $\bar{x}(0)=x_A$ and $\bar{x}(T)=x_B$. It is clear that in general $w_0$ is not equal to $x_A$. 
Let $\bar{t}>0$ and in system \eqref{general_heis_hyst} choose a control $u_1$ such that $(x(\bar{t}),\mathcal{P}[x,w_0](\bar{t}))=(x_A+sign(\dot{\bar{x}}(0))\rho,x_A)$ 
and $u_2(t)\equiv 0$ in $[0,\bar{t}]$. In this way the couple $(x(\bar{t}),\mathcal{P}[x,w_0](\bar{t}))$ will be exactly on one of the two lines $x\pm\rho$ see Figure \ref{initial_data}

Now the value of the input of the play operator and of its output at time $\bar{t}$ are exactly the ones required by Lemma \ref{density}, so for $t>\bar{t}$, let $v^j(t-\bar{t})$ be the function introduced in the preceding Lemma \ref{density}. 
Moreover let $u_1^{j'}(t)$ be the control that steers $v^j(T)$ to $x_B$ in time $\Delta T$, and define the time $T^*=\bar{t}+T+\Delta T$. Then we choose as controls for the hysteretic system \eqref{general_heis_hyst}
\begin{align}
\label{control1}
&\begin{cases}
&u^j_1(t)=u_1(t)\quad\text{for }0\leq t\leq \bar{t}\\
&u^j_1(t)=\frac{\partial v^j}{\partial t} (t-\bar{t})\quad\text{for }\bar{t}\leq t\leq T+\bar{t}\\
&u^j_1(t)=u^{j'}_1(t-T)\quad \text{for } T+\bar{t}<t\leq T+\bar{t}+\Delta T=T^*
\end{cases}\\
\label{control2}
&\begin{cases}
&u^j_2(t)\equiv 0\quad\text{for }0\leq t\leq \bar{t}\\
&u^j_2(t)\equiv \bar{u}_2(t-\bar{t})\quad\text{for }\bar{t}\leq t\leq T+\bar{t}\\
&u^j_2(t)\equiv 0\quad\text{for } T+\bar{t}<t\leq T+\bar{t}+\Delta T=T^*
\end{cases} 
\end{align}
Note that in each of the time intervals $[0,\bar{t}]$ and $[T+\bar{t},T+\bar{t}+\Delta T]$ the variable $x$ is affine and its end-points do not depend on $j$ therefore the two arcs of the trajectories
remain in a compact set. Moreover in the time interval $[\bar{t},T+\bar{t}]$ we can follow the argument sketched in Theorem \ref{approx_control} in the previous section according to which since $\bar{u}_2$ is bounded by hypothesis, $f$ is $C^{\infty}$ and since $v^j$ is bounded (thus also $f(\mathcal{P}[v^j,w_0])$) the hysteretic trajectories $z^j(t)$ do not exit from a common compact set. Therefore in this case also we can consider $f$ bounded and Lipschitz along trajectories.

With the controls chosen in \eqref{control1}-\eqref{control2} we have, according to the definition of $u^{j'}_1$, $x(T^*)=x_B$.
Moreover $y(T^*)=y_B$ since we are using exactly the control $\bar{u}_2$ up to time $T$ and then zero.
Finally 
\begin{equation}
\small
\begin{aligned}
&|z(T^*)-\bar{z}(T)|=	\\& \Bigl |\int_{0}^{\bar{t}} 0\,dt+\int_{\bar{t}}^{T+\bar{t}}\bigl(f(\mathcal{P}[v^j,w_0](t-\bar{t}))-f(\bar{x}(t-\bar{t}))\bigr)\bar{u}_2(t-\bar{t})\,dt+\int_{T+\bar{t}}^{T^{*}}0\, dt\Bigr |\\
&\leq \int_0^T\bigl | \bigl(f(\mathcal{P}[v^j,w_0](t-\bar{t}))-f(\bar{x}(t-\bar{t}))\bigr)\bigr | \|\bar{u}_2\|_{\infty}\,dt\\
&\leq L T\|\bar{u}_2\|_{\infty} \|\mathcal{P}[v^j,w_0]-\bar{x}\|_{\infty}\to 0 \,\,\text{as $j\to\infty$}
\end{aligned}
\end{equation}
where $L$ is the Lipschitz constant of the function $f$ and we have used the convergence of $\mathcal{P}[v^j,w_0]$ to $\bar{x}$ in $L^{\infty}$ given by Lemma \ref{density}.
\EndProof
\begin{proposition}
In the case in which $f(x)=x$, i.e. the classical Heisemberg system, we have the exact controllability in the hysteretic case.
\end{proposition}
\Proof
First of all note that, according to \eqref{lie_heis} in order to move only in the $z$ direction the strategy is to move along the direction of the Lie bracket of the two vector fields. It is well known that this can be achieved choosing for example controls 
\begin{equation}
\label{control_loop}
\begin{aligned}
u_1=\begin{cases}
\alpha\quad 0\leq t\leq T\\
0\quad T\leq t\leq 2T\\
-\alpha\quad 2 T\leq t\leq 3T\\
0\quad 3T\leq t\leq 4T
\end{cases}&&u_2=\begin{cases}
0\quad 0\leq t\leq T\\
\beta\quad 2 T\leq t\leq 3T\\
0\quad 2 T\leq t\leq 3T\\
-\beta\quad 3T\leq t\leq 4T
\end{cases}.
\end{aligned}
\end{equation} 
and that, after a period of $4T$, we have (see, for example, Coron \cite{Coron56}, page 130) \begin{equation}
\label{Lie_bracket}
\Delta z=z(4T)-z_A=T^2\alpha\beta
\end{equation}

\noindent
This means that, starting from $(x_A,y_A,z_A)$ and using the controls in (\ref{control_loop}), after the time interval $4T$ the system is in position $(x_A,y_A,z_A+T^2\alpha\beta)$.

At first we are interested in studying the controllability properties of the hysteretic version of \eqref{general_heis}
\begin{equation}
\label{heis_hyst}
\begin{pmatrix}
\dot x\\
\dot y\\
\dot z
\end{pmatrix}=\begin{pmatrix}1\\0\\0\end{pmatrix}u_1+\begin{pmatrix}0\\1\\\mathcal{P}[x,w_0]\end{pmatrix}u_2
\end{equation}
where $\mathcal{P}[x,w_0]$ denotes the play operator applied to the real variable $x$ with initial datum $w_0$. We will show that system \eqref{heis_hyst} is controllable.
The proof is based on the surjectivity of the play operator on a certain class of functions. 
In a first time interval $[0,t_1]$ we use a control $u_2$ which drives the variable $y$ from $y_A$ to $y_B$ and $z$ from $z_A$ to a certain $\bar{z}$ and we set $u_1=0$. Now observe that, for suitable $\tilde\alpha,\tilde\beta$, we can consider two controls as in \eqref{control_loop}, which we call $\tilde{u}_1,\,\tilde{u}_2$, such that in an interval of time $\Delta t=4T$ drive the non hysteretic system \eqref{general_heis} from $\bar{z}$ to $z_B$ and drive $x$ and $y$ back to their initial values. Also note that, assuming $\overline z\neq z_B$ (otherwise the variable $z$ is already setted), we certainly have $\tilde\alpha\neq0$.
For the hysteretic Heisenberg system \eqref{heis_hyst}, we have first to suitably arrange the initial datum: in a time interval $[t_1,t_2]$ we switch off $u_2$ and turn on $u_1$ in such a way that the pair $(x(t_2),w(t_2))=(x(t_2),\mathcal{P}[x,w_0](t_2))$ is in one of the two boundary lines of $\Omega_\rho$, in particular $w(t_2)=x(t_2)-\rho$ if $\tilde{u}_1(0)>0$ or $w(t_2)=x(t_2)+\rho$ if $\tilde{u}_1(0)<0$ (see Figure \ref{initial_data}).
Now, note that if we translate the controls $\tilde{u}_1,\,\tilde{u}_2$ to the time interval $[t_2,t_2+4T]$ and integrate them, we obtain continuous piecewise linear functions which belong to a particular class of continuous functions on which we have already proved that the play operator is surjective (see Remark \ref{surjectivity} regarding Lemma \ref{approx_controls}). Therefore, denoting $\tilde x(t)=x(t_2)+\int_{t_2}^t\tilde u_1(s)ds$ for all $t\in[t_2,t_2+4T]$, 
it is possible to find a control $v_1$ in $[t_2,t_2+4T]$ such that the output of the play operator acting on the corresponding trajectory $x$ ($\dot x=v_1$) gives exactly $\tilde{x}$, i.e 
$$
\exists v_1\quad s.t.\quad\ \forall\ t\in[t_2,t_2+4T],\ x(t)=\int_{t_2}^{t} v_1(s)\,ds\quad \Rightarrow \mathcal{P}[x,w(t_2)](t)=\tilde{x}(t)
$$
Now using this control $v_1$ and the control $\tilde{u}_2$ corresponding to the trajectory $t\mapsto\tilde{y}=y_B+\int_{t_2}^t\tilde u_2(s)ds$ in the time interval $[t_2,t_2+4T]$, the system moves from $(x_A,y_B,\bar{z})$ to $(x(t_2+4T),y_B,z_B)$, . Finally it suffices to set $u_2=0$ and use a control $u_1$ that adjusts the $x$ variable in a time interval $[t_2+4T,t_3]$ in order to get $x(t_3)=x_B$.
\EndProof
\begin{remark}
Note that the strategy proposed for system \eqref{heis} is not valid for a generic Lipschitz function $f$, indeed for a generic nonlinear $f$, a control loop like \eqref{control_loop} does not lead to a displacement \eqref{Lie_bracket} for any $T$ but only for small ones. Thus it is only a local approximation of the displacement for small times, and cannot be used to prove the exact controllability of the system, for any initial and final position. Nevertheless the approximating strategy is still valid and can be used to get at least the approximate controllability result.
\end{remark}
\begin{remark}
\label{rmrk:previous}
One can think to use directly the limit of the controls $\vecu^j$ for the hysteretic system, but the sequence of controls $u_1^j$ converges only point-wise to the piecewise constant control $\bar{u}_1$, and this control is not good for the hysteretic system. Indeed, using it, we are not able to reproduce the trajectory $\bar{x}$ as the output of the play operator. Indeed, the play operator is continuous only with respect to the topology of uniform convergence but we can grant only pointwise convergence.\\
Moreover note that the result obtained using the extended definition of the Play operator for discontinuous inputs (\cite{brospr,Kreici,VR1,VR2}), would require the use of measure control instead of measurable ones. More precisely the good discontinuous input $v$ such that $P[v,w_0]=\bar x$, can be obtained using a control which is $\bar u_1+2\rho\sum_{\text{number of jumps}}\delta_{t_i}$, where $\delta_{t_i}$ is the Dirac delta function centered in the jumping times. This control is a measure which is not the class of controls in which we are more interested, that are the measurable controls.
\end{remark}
\CorrR{
\begin{remark}
\label{rmrk:bound}
Note that the controls in \eqref{control1} \CorrB{are not equi-bounded} in $L^{\infty}$ as $j\to\infty$ but only in $L^1$. However this is not a problem for our construction since our control generates the input of the play operator and we do not require that it converges strongly but only that its output does. Moreover our strategy shows how to move, in the limit $j\to\infty$, between two fixed points $A$ and $B$ not how to reproduce the non-hysteretic trajectory $\bar{x}$, compare with Remark \ref{rmrk:previous}. Instead using $u_1^j$ we generate a trajectory $v^j$ whose output through the play operator strongly converges to $\bar{x}$.\\
Moreover, note that if the admissible controls must take value in a compact set (i.e. they are equibounded), then our approximation procedure, even if we use the extended version of the play operator to discontinuous inputs, is not more applicable, since the controls $u_1^j$ are clearly unbounded as $j\to\infty$. Indeed, with controls in a compact set we only have the estimate  
$$
||\mathcal{P}[v^j,w_0]-\bar{x}||_\infty\geq \frac{1}{M}(\max_i|\alpha_i|),
$$

\noindent
where $M$ is the bound on the controls (see also the estimate just under Figure \ref{traj_approx}). Hence, even in the limit, only a partial controllability result may in general hold.
\end{remark}
}
\begin{remark}
Note that  as in remark \ref{Robustness} also in this case we have a sort of robustness of the system. 
\end{remark}
\subsection{Generalization}

In the previous paragraph we have considered "triangular" systems of the form \eqref{general_heis_hyst}, and that structure was crucial in the proof of the result. The controllability of different systems with hysteresis in the space is still under investigation. However, here we give a controllability result for a further generalization of that particular "triangular" structure, to the case of more hysteresis dependent variables. More precisely, we achieve it  iterating the procedure previously described.

Let us start from the following system
\begin{equation}
\begin{pmatrix}\dot{x}_1\\\dot{x}_2\\\dot{x}_3\\\dot{y}_4\\\dot{y}_5\end{pmatrix}=\begin{pmatrix}1\\0\\0\\0\\0\end{pmatrix}u_1+\begin{pmatrix}0\\1\\0\\f_2(\mathcal{P}[x_1,w_0^1])\\0\end{pmatrix}u_2+\begin{pmatrix}0\\0\\1\\0\\f_3(\mathcal{P}[x_1,w_0^1],\mathcal{P}[x_2,w_0^2])\end{pmatrix}u_3
\end{equation}
Suppose that $\bar{\vecu}=(\bar{u}_1,\bar{u}_2,\bar{u}_3)$ is the piecewise constant control that steers the non hysteretic system from $A$ to $B$ and let $\bar{x}_1(t),\bar{x}_2(t)$ be the corresponding continuous piecewise linear trajectories. The idea is to mimic what we have done in the previous example. Therefore the first thing to do is to bring $(x_1,w^1)$ and $(x_2,w^2)$ in the good position, like in Figure \ref{initial_data}, i.e. the points $(x_1,w^1)$ and $(x_2,w^2)$ in the graph $(x,w)$ have to be on the boundary of the set $ \Omega_\rho$ (see \eqref{omega_rho}). Thus we put $u_3\equiv 0$ and we choose suitable $u_1$ and $u_2$ such that  $(x_1(\bar{t}),w^1(\bar{t}))=(x_{1_A}+sgn(\dot{\bar{x}}_{1_A})\rho,x_{1_A})$ and $(x_2(\bar{t}),w^2(\bar{t}))=(x_{2_A}+sgn(\dot{\bar{x}}_{2_A})\rho,x_{2_A})$. Then using Lemma \ref{density} we reproduce the trajectory $\bar{x}_1(t),\bar{x}_2(t)$ using a sequence of controls $(u_1^j,u_2^j)$ that produce the play input functions $(v_1^j,v_2^j)$. In this way at a certain time $\bar{t}+T$ we have that $x_3(\bar{t}+T)=x_{3_B}$ and $y_5(\bar{t}+T)\to y_{5_B}$.
More precisely

\begin{equation}
\small
\begin{array}{ll}
&u^j_1(t)=u_1(t)\\\\
&u^j_2(t)=u_2(t) \\\\
&u^j_3(t)\equiv 0
\end{array} \text{for }0\leq t\leq \bar{t}
\begin{array}{ll}
&u^j_1(t)=\frac{\partial v^j_1}{\partial t} (t-\bar{t})\\\\
&u^j_2(t)=\frac{\partial v^j_2}{\partial t} (t-\bar{t})\\\\
&u^j_3(t)\equiv \bar{u}_3(t-\bar{t})
\end{array} \text{for }\bar{t}\leq t\leq T+\bar{t}
\end{equation}

 Now it is important to note that for the structure of the vector fields the direction $\partial_{x_1}$, $\partial_{x_2}$ and $\partial_{y_4}$ can be generated using only the first two controls. Therefore we now consider only the subsystem
\begin{equation}
\begin{pmatrix}\dot{x_1}\\\dot{x_2}\\\dot{y_4}\end{pmatrix}=\begin{pmatrix}1\\0\\0\end{pmatrix}u_1+\begin{pmatrix}0\\1\\f_2(\mathcal{P}[x_1,w_1(\bar{t}+T)]\end{pmatrix}u_2
\end{equation}
This system is exactly the one of the example for which we already showed that it is possible to find a sequence of controls that steers it from $(x_1(\bar{t}+T),x_2(\bar{t}+T),y_4(\bar{t}+T))$ to $(x_{1_B},x_{2_B},y_{4_B})$.

It is possible to generalize this idea and consequently the proof of Theorem \ref{thm3} to systems with the following coordinates:  $\vecz:=(x_1,\cdots,x_m,y_{m+1},\cdots,y_{2m-1})$ and control vector fields $\vecg_i$ of the following type
\begin{equation*}
\begin{aligned}
&\vecg_1(\vecz):=\partial_{x_1}\\
&\vecg_i(\vecz):=\partial_{x_i}+f_{i}(x_1,\cdots,x_{i-1})\partial_{y_{m+i-1}}\quad i=2\cdots m
\end{aligned}
\end{equation*}
Suppose all $f_i\in C^{\infty}$, and that the vector fields are bracket generating and thus the associated control system is controllable.

The hysteretic system that we want to investigate is now
\begin{equation}
\dot{\vecz}=\sum_{i=1}^m\vecg_i(\mathcal{P}[\vecz,\vecw_0])u_i
\end{equation}
where by $\mathcal{P}[\vecz,\vecw_0]$ we mean the component wise scalar play operator.
To move the system between two fixed points $A:=\vecz_A$ and $B:=\vecz_B$ it suffices to iterate the preceding procedure.
The idea is the following: suppose that $\bar{\vecu}$ is the piecewise constant control that steers the non hysteretic system from $A$ to $B$ and let $\bar{x}_1(t),\cdots,\bar{x}_{m-1}(t)$ the corresponding trajectory of the first $m-1$ coordinates. For the structure of the vector fields these $\bar{x}_i(t)$ are continuous piecewise linear functions. After reaching the input-output relation between $(x_1,\cdots, x_{m-1})$ 
and $(w_1,\cdots,w_{m-1})$ given by Figure \ref{initial_data} (i.e. $(x_i,w_i)\in\partial\Omega_\rho$ for all $i=1,\dots,m-1$), since we know that the play operator has dense image on the continuous piecewise linear functions, we are able to find sequences $v^j_i(t)$ such that their output $\mathcal{P}[ v^j_i](t)$ converge to $\bar{x}_i(t)$ as $j\to\infty$, and use the corresponding controls $u^j_i(t):=\frac{\partial v^j_i}{\partial t}(t)$ for $i=1,\dots m-1$ and $u^j_m(t):=\bar{u}_m(t)$. In this way for sure $y_{2m-1}\to y_{2m-1}^B$, but the other coordinates can be different. Thus the strategy now is to put $u_m\equiv 0$ and find controls $u_1,\cdots u_{m-1}$ which adjust the coordinate $y_{2m-2}$ in the non hysteretic system. This is possible since the chosen
vector fields structure allows to generate the $y_{m+i-1}$-direction using only the first $i$ controls.
\begin{remark}
Observe that even if the type of vector fields for which the preceeding procedure works seems to be restricted, they belong to Carnot groups of step $2$. These groups are widely used in sub-riemannian geometry and control theory (see Agrachev et al \cite{Agrachev12}), and describe a wide class of mechanical systems, starting from the Heisenberg flywheel system to its other  generalization.
\end{remark}
\begin{remark}
Note that the essential property of the Play operator used to build the approximating sequence of controls is Lemma \ref{density}, i.e. the fact that it has dense image in the space of continuous piecewise linear functions. This means that the theorem on the approximate controllability is valid also tanking other hysteresis operators. An example is the so called sweeping process (see Moreau and Colombo et al. \cite{Colombo, Moreau}) which is built as follows. Consider a moving set $C(t)$ in $\mathbb{R}^n$, depending on the time $t\in[0,T]$, and an initial condition $\vecz_0\in C(0)$. In several contexts, the modelization of the displacement $\vecz(t)$ of the initial condition $\vecz_0$ subject to the dragging, or sweeping due to the displacement of $C(t)$ pops up. It is natural to think that the point $\vecz(t)$ remains at rest until it is caught by the boundary of $C(t)$ and then its velocity is normal to $\partial C(t)$. It is a kind of one sided movement. Formally, the sweeping process is the differential inclusion with initial condition
$$\vecz(t)\in-N_{C(t)}(\vecz(t)),\qquad\vecz(0) =\vecz_0\in C(0)$$
where $N_{C}(\vecz)$ denotes the normal cone to $C$ at $\vecz\in C$. The sweeping process has the property of having  dense image in the space of continuous piecewise linear functions. This is clear since in dimension one it behaves exactly as the play operator  (see Recupero et al. \cite{Recupero}), more precisely it is the case when $C(t)$ is a translation of a symmetric closed convex set $C$:
$$
\dot z(t)\in-N_{u(t)+C}(z(t))
$$
\end{remark}
\subsection{The case of switching hysteresis}
Again, we consider the system
\[
\dot \vecz=\sum_{i=1}^m\vecg_i(\vecz)u_i,
\]
but now we suppose that every $\vecg_i$ incurs in a discontinuity across an hyperplane of $\mathbb{R}^n$. In particular (because of an approximation point of view as well as a possible intrinsic hysteretic behavior, see Liberzon \cite{lib}) we describe such a discontinuity by a delayed relay.
For every $i=1,\dots,m$ let $\xi_i\in\mathbb{R}^n$ be a unit vector, representing the unit normal to the discontinuity hyperplane, and we consider the delayed relay with hysteresis as in the figure.
\begin{figure}[H]
\centerline{
\includegraphics[scale=0.5]{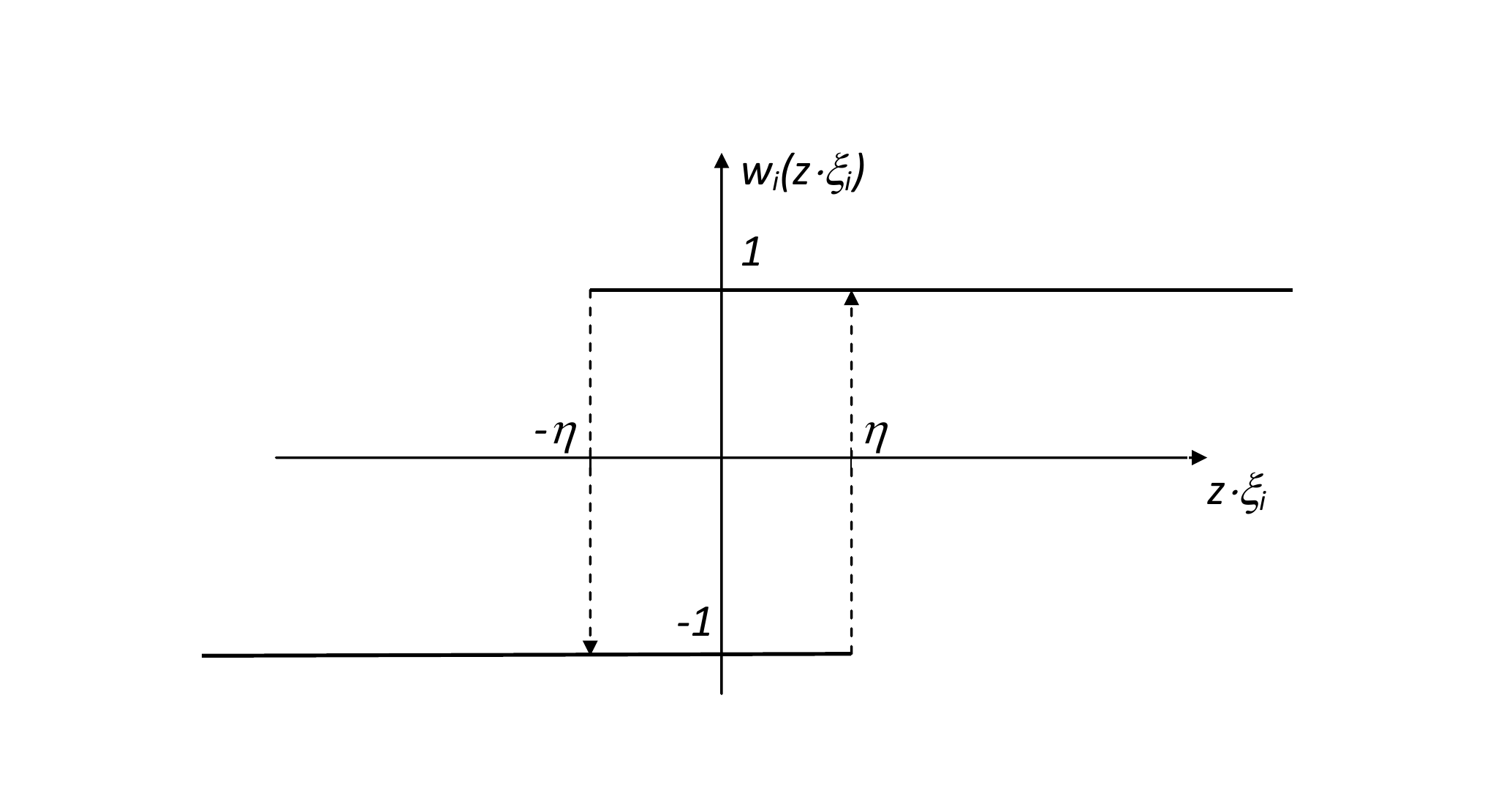}
}
\end{figure}
For every $i=1,\dots,m$ and for every $w_i\in\{-1,1\}$, we have a field $\vecg_i^{w_i}$. We then consider the controlled systems
\begin{equation}
\label{switching_dyn}
\dot \vecz=\sum_{i=1}^m\vecg_i^{w_i(\vecz\cdot\xi_i)}(\vecz)u_i,
\end{equation}
that is, each field $\vecg_i$ is subject to switch, in dependence on $\vecz\cdot\xi_i$, with a delayed rule. 
\begin{hypo}
\label{LARC}
 For every $m$-string $(w_1,w_2,\dots,w_m)\in\{-1,1\}^m$, the $m$ fields $(\vecg_1^{w_1},\vecg_2^{w_2},\dots,\vecg_m^{w_m})$ satisfy Assumption \ref{hypo_controllability}, i.e. their Lie algebra is fully generated.
 \end{hypo}
 The controllability question is now whether, given, $A,B\in\mathbb{R}^n$, there are or not piece-wise constant controls $u_i$ which steer the system from $A$ to $B$ in a finite time $T>0$.

\begin{figure}[H]
\centering
\includegraphics[scale=0.45]{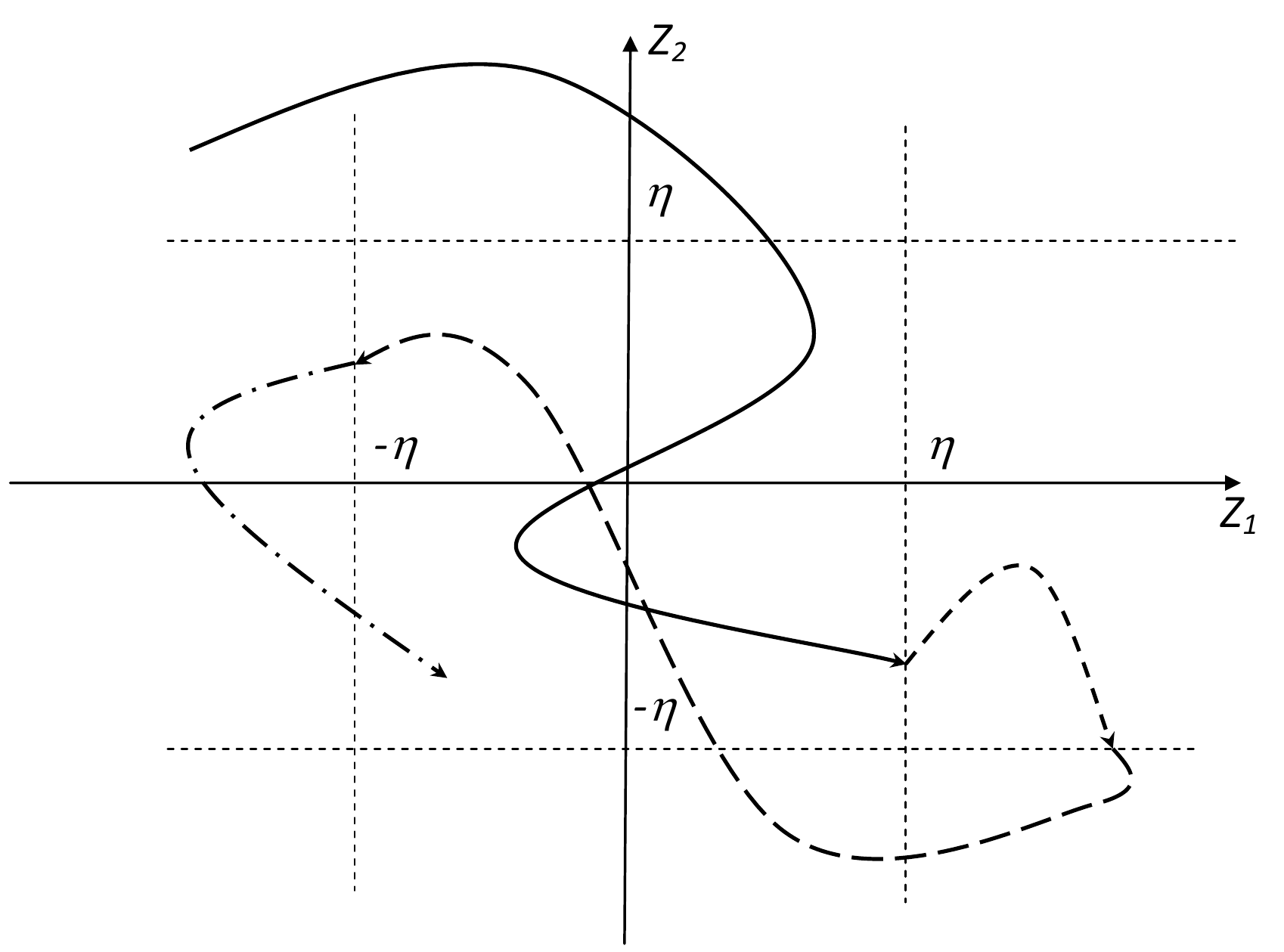}
\caption{\label{projection}The projection of the trajectory on the first two coordinates}
\end{figure}

Just as an example, suppose that $m=2$ and that $\xi_1$ and $\xi_2$ are respectively the first and second element of the canonical basis of $\mathbb{R}^n$. Then we may have the situation of Figure \ref{projection}, for the projection of the trajectory on the first two coordinates, where, for a given control $\vecu=(u_1,u_2)$,  the filled curve is the evolution with $(\vecg_1^{-1},\vecg_2^1)$, the short dashed curve is the evolution with $(\vecg_1^1,\vecg_2^1)$, the long dashed curve is the evolution with $(\vecg_1^1,\vecg_2^{-1})$ and the point-dashed one is the evolution with $(\vecg_1^{-1},\vecg_2^{-1})$.
The state space $\mathbb{R}^n$ is then divided in $2^m$ (non-disjointed) sectors, every one indexed by the corresponding $m$-string of $1$ and $-1$. For example, with respect to the figure, we have the sector indexed by $(1,-1)$ which is $[-\eta,+\infty[\times]-\infty,\eta]\times\mathbb{R}^{n-2}$. 
When we start to move inside one of the sectors, then we continue to move in the same mode $(\vecg_i^{w_i})_i$ until we leave that sector, and after that we move in the new modality (corresponding to the index of the new sector) determined by the delayed switching rule. Since the sectors have non-empty intersection, then, together with the starting point $z_A\in\mathbb{R}^n$, whenever it belongs to more than one sector, we must also give the initial sector (i.e. the initial index $(w_i)_i$, i.e. the initial evolution mode).
Note that, every point $z_B\in\mathbb{R}^n$ stays in the interior of a sector and hence, in our controllability problem, we can always equip $z_B$ with the index of that sector. More precisely, in the following, we are going to decide to reach $z_B$ with exactly that mode of evolution in the last part of the time interval.
\begin{remark}
If the discontinuity is not ``delayed", then the trajectory may not even exist when using piecewise constant controls. Indeed, consider the following simple example. Suppose $m=2$ $\xi_i=e_i$, $\vecg_i^1=e_i=-\vecg_i^{-1}$. Now, for each $w_i$, we have the unique switching threshold $z_i=0$ (not delayed). Hence the sectors (after projection on the first two coordinates) are the quadrants of $\mathbb{R}^2$ with intersections on the axes. Then, for example, starting from the sector $(1,1)$, it  is impossible to leave it just using piecewise controls. Indeed, 
take for example (the projection of) $A=(0,1)$ and $(w_1,w_2)=(1,1)$ (the first quadrant). To leave, we have to take $u_1=-1$ because $\vecg_1^1=(1,0)$. But then immediately $\vecg_1^{-1}$ will become $(-1,0)$ and hence we have no existence of the trajectory if $u_1$ is constant in at least small time interval $[0,\delta]$. Anyway note that, even in the case of not delayed discontinuity we may still have controllability. But we may be forced to consider some state-dependent restrictions of the set of admissible controls, and then the controllability conditions may be less immediate and useful. Instead, considering a ``delayed discontinuity'' allows us to still use all the set of admissible controls and then to obtain a more standard and manageable set of controllability conditions. 
\end{remark}
\begin{theorem}
\label{thm:switching}
Given Assumption \ref{LARC}, system \eqref{switching_dyn} is controllable. That is, for every $A,B\in\mathbb{R}^n$, and for every initial $m$-string $(w_1,\dots,w_m)\in\{-1,1\}^m$, compatible with $A$ (i.e. $A$ belongs to the corresponding sector), there exist $T>0$ and piece-wise constant controls $u_i$ which steer the system from $A$ to $B$ in time $T$.
\end{theorem}
\Proof
As already said, we can always consider $B$ as belonging to the interior of one sector and decide to reach it with that mode of evolution. Note that the initial sector of $A$ is given by the initial values of $w_i$ which are not at our disposal.  We distinguish various cases depending on $A$.

{\it First case}: $A$ and $B$ belong to the interior of same sector $S$. We look for an admissible trajectory connecting them without leaving $S$. Since each sector is connected we can follow the construction proposed in Laumond et al. \cite{Laumond}. At first we consider a continuous path $\gamma:[0,1]\to\mathbb{R}^n$ of finite length (not necessarily admissible), connecting $A$ to $B$ and not leaving the interior of the sector. Let us define $\delta:=\min_{t\in[0,1]} dist(\gamma(t),\partial S)>0$. Let $K\subset \mathbb{R}^m$ be a closed ball centered at the origin and, for every $\vecz\in\mathbb{R}^n$ and $T>0$, define $R_\vecz (T)$ the set of configurations reachable from $\vecz$ by an admissible trajectory before the time $T$ and only using $\vecu\in K$. For any $\vecz\in\gamma([0,1])$ there exists $T_\vecz$ such that $R_\vecz(T_\vecz)$ is contained in $S$ and since each $R_\vecz(T_\vecz)$ is open \CorrR{(by the Chow condition Assumption \ref{hypo_controllability}, see Coron \cite{Coron56} )} it contains a ball of radius $\epsilon_\vecz$, $\mathcal{B}(\vecz,\epsilon_\vecz)$, with $\epsilon_\vecz<\delta$ so that it does not intersect $\partial S$. Let us take a finite covering by $N$ balls of this kind of the compact set $\gamma$, whose centers 
belongs to $\gamma$.  We may suppose that $\vecz_1=A=$ and $\vecz_N=B$ are in such set of centers. 
We can also arrange the labels in such a way that $\vecz_2\neq\vecz_1$ and $\mathcal{B}(\vecz_1,\epsilon_{\vecz_1})\cap \mathcal{B}(\vecz_2,\epsilon_{\vecz_2})\neq \emptyset$, otherwise we would have not a covering of $\gamma$ with open balls centered in $\gamma$.
Similarly, if $\vecz_2\neq B$ we may arrange that $\vecz_3\neq \vecz_2,\,\vecz_1$ and $\mathcal{B}(\vecz_2,\epsilon_{\vecz_2})\cap \mathcal{B}(\vecz_3,\epsilon_{\vecz_3})\neq \emptyset$. We proceed in this way. For any $i$ take $\vecz_{i,i+1}\in\mathcal{B}(\vecz_i, \epsilon_{\vecz_i}) \cap \mathcal{B}(\vecz_{i+1}, \epsilon_{\vecz_{i+1}})$.  Thus there is an admissible trajectory from $\vecz_i$ to $\vecz_{i,i+1}$ and another one from $\vecz_{i,i+1}$ to $\vecz_{i+1}$ that do not exit from $R_{\vecz_i} (T_{\vecz_i} ) \cup R_{\vecz_{i+1}} (T_{\vecz_{i+1}} )\subset S$. The sequence $(\vecz_i)_{i=1}^N$ is finite and we can conclude that there exists an admissible path from $A$ to $B$, running for a time $T\le\sum_iT_{\vecz_{i}}$, that does not exit the sector and, in particular, that allows us to not switch dynamics.

{\it Second case:} $A$ and $B$ do not belong to the same sector, but $A$ is internal to its initial sector. We restrict to the case where the sectors of $A$ and $B$ differ by one switching only (their switching $m$-dimensional labels differ by one component only). The other cases can be constructed iteratively in a similar way. In this case we make a first step: starting from $A$ we use a piece-wise constant control which allows us to reach the switching boundary between the sector of $A$ and the one of $B$ and without leaving the sector of $A$. In particular, we may reach a point which is not a ''multiple-switching'' point, but a ''single-switching'' point between the sector of $A$ and $B$ (referring to the example in Figure  \ref{projection}, it is not a corner point of the kind $(\pm\eta,\pm\eta)$ where both switching may occur).
This can be done following the strategy of step one until an interior point sufficiently near the boundary, and then, again thank to the Chow condition, we can use a dynamics that allows us to reach the boundary and then switch. We are now in the situation of the previous case, since the point on the boundary on which we are arrived, after the switching is now an internal point of the sector of $B$. Note that without the delayed relay we do not switch in the interior of the other sector. 

{\it Third case:} $A$ belongs to the boundary of its initial sector. Again, thank to the Chow condition (Assumption \ref{LARC}), we may initially use a control that makes the trajectory switch and so almost immediately reach a point  in the interior of one of the sectors. Then we proceed as in one of the previous cases. 
\EndProof
\subsection{Switching and play hysteresis: approximation and controllability}
\label{S_P_Hysteresis}
In the previous subsection we have treated the case where the system is affected by a delayed switching hysteresis, with discontinuous output. Let us note that the play hysteresis is instead a continuous hysteresis (the output is continuous). Actually, the play operator can be seen as a superposition of an infinitely many quantity of delayed relays, and then can be in some sense approximated by a big, but finite, number of delayed relays.

For any $r\in[0,1]$ let $h_r$ represent the relay with threshold $(-1+r,r)$. For a given scalar time-continuous input $\zeta$, and for given initial output states for each relays (which we do not display for simplicity of notations) let us consider the following "macroscopic" output

\[
w(t)=\int_0^1h_r[\zeta](t)dr.
\]

\begin{figure}[H]
\centerline{
\includegraphics[scale=0.6]{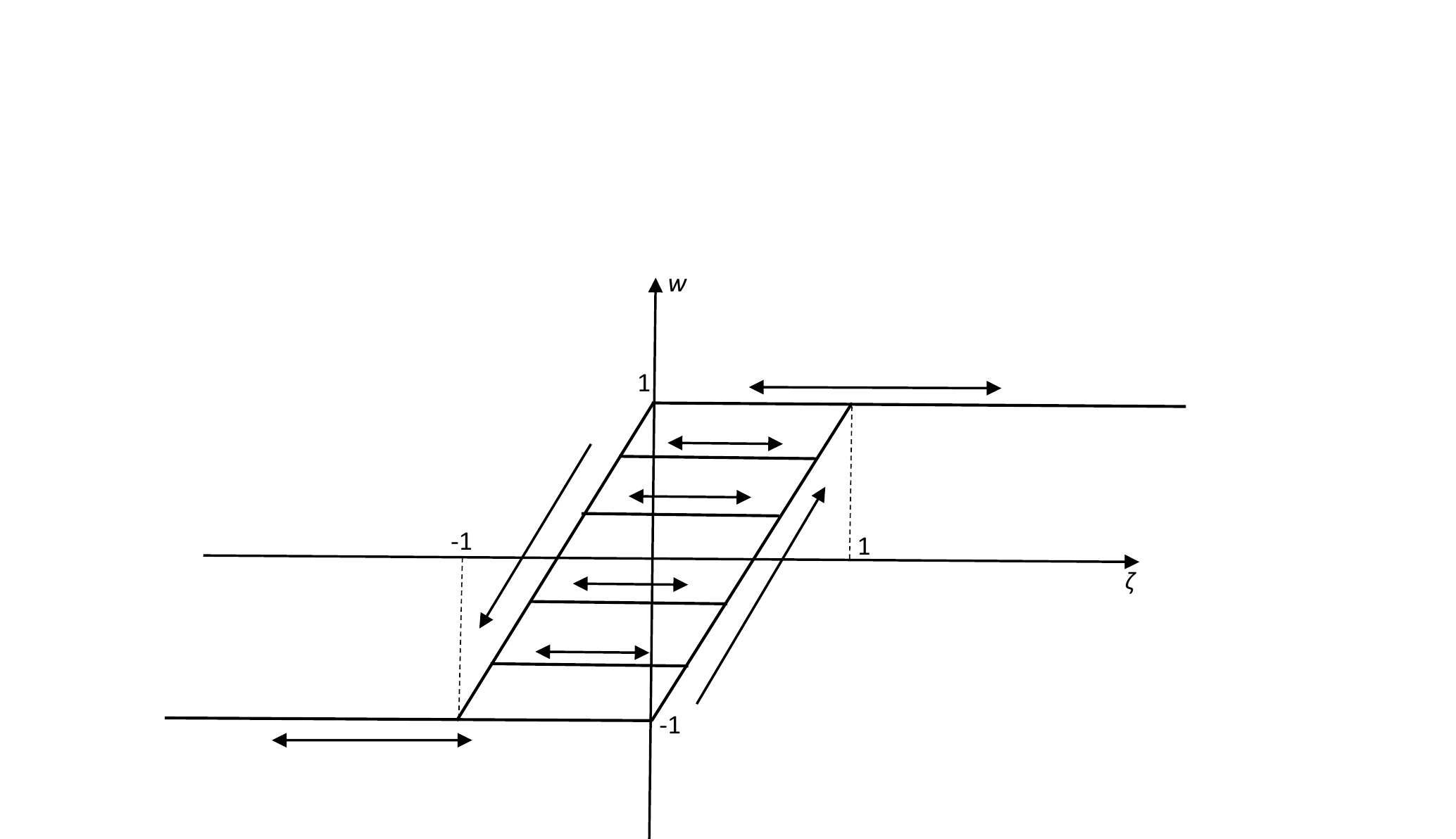}
}
\caption{\label{fig:truncated_play} Truncated play hysteresis-loop}
\end{figure}

\noindent
Since any relays is identified by $r\in[0,1]$, the output states of the relays is, at any time $t$, a function $o^t:[0,1]\to\{-1,1\}$, $r\mapsto h_r[\zeta](t)$. We consider the following hypothesis 
\begin{equation}
\label{eq:hysteron}
\exists\ r^\tau\in[0,1]\ \mbox{such that } o^\tau(r)=1\ \mbox{if } r<r^\tau,\  o^\tau(r)=-1\ \mbox{if }r>r^\tau.
\end{equation}
\noindent
If (\ref{eq:hysteron}) is satisfied for some $\tau\ge0$, then, subject to the evolution of $\zeta$, the pair input-output $(\zeta,w)$ evolves, for $t\ge\tau$, inside the following hysteresis-loop, with the described evolution  by the arrows, which exactly corresponds to a truncated Play operator with slope $2$ and width $\rho=1$ (see Figure \ref{fig:truncated_play}).

\noindent
In particular, if (\ref{eq:hysteron}) holds, then for every time $t\ge\tau$, there exists $r^t\in[0,1]$ such that $o^t(r)=1$ for $r<r^t$ and $o^t(r)=-1$ for $r>r^t$, i. e. (\ref{eq:hysteron}) holds. If instead the initial output does not satisfy (\ref{eq:hysteron}), then the evolution of $(\zeta,w)$ is not necessarily described by the hysteresis-loop as in Figure \ref{fig:truncated_play}, but, whenever at a time $t$, the hypothesis is satisfied, then the evolution will remain inside that hysteresis loop for all subsequent times. Also note that, if at a certain time $t$ it is $\zeta(t)\ge1$ (respectively, $\zeta(t)\le-1$), then all the relays are switched on $1$ (respectively, $-1$) and (\ref{eq:hysteron}) is satisfied. Hence, acting if necessary on the input $\zeta$, we can always suppose to start the evolution satisfying (\ref{eq:hysteron}). 

Now, instead of considering a continuum of relays indexed by $r\in[0,1]$, we consider $k$ relays, $k\in\mathbb{N}\setminus\{0\}$, $h_1,\dots,h_k$, with the hypothesis that, for every $i=1,\dots,k$, $h_i$ has thresholds $(-1+i/k,i/k)$. We then consider the macroscopic output

\begin{equation}
\label{eq:w_n}
w_k[\zeta](t)=\frac{1}{k}\sum_{i=1}^kh_i[\zeta](t).
\end{equation}

\begin{figure}[H]
\centerline{
\includegraphics[scale=0.6]{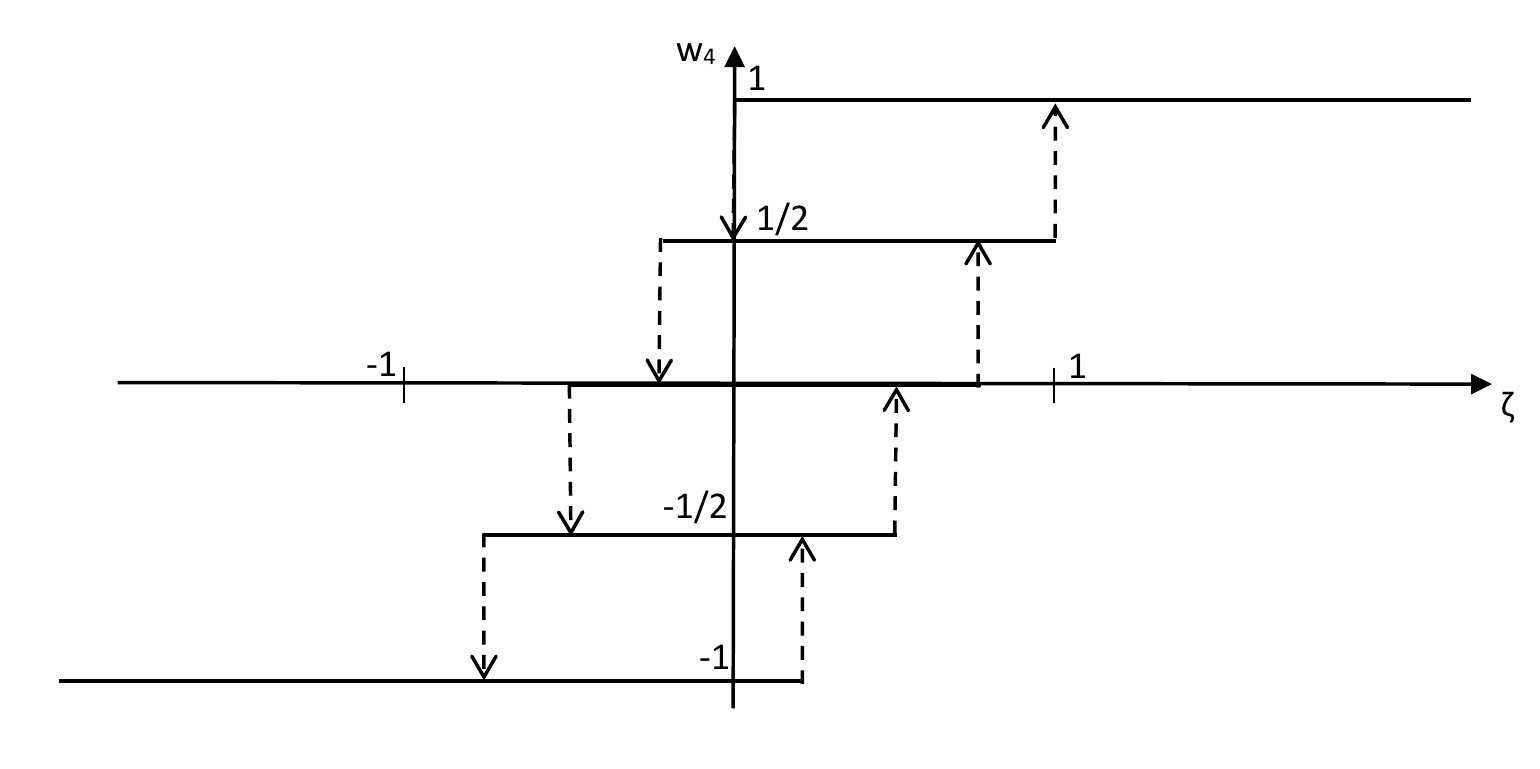}
}
\caption{\label{fig:play_discreto} Here we have represented the case of $k=4$ relays, with thresholds $(-3/4,1/4), (-1/2,1/2), (-1/4, 3/4), (0,1)$.}
\end{figure}

Now, at any time $t$, the output states of the $k$ relays is a function $o_k^t:\{1,\dots,k\}\to\{-1,1\}$. Similarly to (\ref{eq:hysteron}), we consider the following hypothesis:
\begin{equation}
\label{eq:hysteron_k}
\exists\ i_k^\tau\in\{1,\dots,k\}\ \mbox{such that } o_k^\tau(i)=1\ \mbox{if } i<i_k^\tau,\ o_k^\tau(i)=-1\ \mbox{if } i>i_k^\tau. 
\end{equation}
\noindent
If (\ref{eq:hysteron_k}) is satisfied for some $\tau\ge0$, then the the pair $(\zeta,w_k)$ evolves, for $t\ge\tau$, inside the following discrete hysteresis-loop (see Figure \ref{fig:play_discreto}) which evidently approximates the continuous hysteresis-loop of Figure \ref{fig:truncated_play}.

\noindent 
In this case, it can be seen that, if we start from an initial output states that does not satisfy (\ref{eq:hysteron_k}), then after a finite number of switches, $(\zeta,w_k)$ necessarily enters the hysteresis-loop of Figure \ref{fig:play_discreto} and will remain there for all the times. Moreover note that, when we are in that hysteresis-loop, at any time $t$ (\ref{eq:hysteron_k}) is satisfied and  there is a bijection between the possible values of the output $w_k$ and the $k$-tuple $(1,\dots,1,-1,\dots,-1)=(h_1[\zeta](t),\dots,h_k[\zeta](t))$ of the relays outputs, image of the function $o_k^t$.

Obviously, the larger $k$, the better the discrete hysteresis loop in Figure \ref{fig:play_discreto} shape, approximates the shape of the continuous one in Figure \ref{fig:truncated_play}.

We now consider the controllability of the system

\begin{equation}
\label{eq:play_truncated_system}
\dot z=\sum_{j=1}^mg_j({\cal P}[z\cdot\xi_j],z)u_j
\end{equation}

\noindent
where $\cal P$ is the truncated play operator as in Figure \ref{fig:truncated_play}, and $\xi_j$, $j=1,\dots,m$ are linearly independent unit vectors.  Moreover, we assume the following hypothesis
\begin{equation}
\small
\begin{aligned}
\label{eq:hypothesis_versus_truncated}
&\forall\ w_j\in[-1,1]\ \mbox{the $m$ fields } z\mapsto g_j(w_j,z) \mbox{ satisfy the the controllability}\\
&\mbox{ hypothesis, i.e. their Lie algebra is fully generated (Assumption \ref{hypo_controllability})} .
\end{aligned}
\end{equation}

\noindent
We take $k\in\mathbb{N}\setminus\{0\}$, and replace the truncated play operator $\cal P$ by $w_k$ as in (\ref{eq:w_n}). We then get the "approximating problem" of controlling the system

\begin{equation}
\label{eq:discrete_system}
\dot z=\sum_{j=1}^mg_j(w_k[z\cdot\xi_j],z)u_j
\end{equation}

\begin{theorem}
\label{thm:multi_switching_controllability}
System (\ref{eq:discrete_system}) is exactly controllable.
\end{theorem}

\begin{figure}[H]
\centerline{
\includegraphics[scale=0.7]{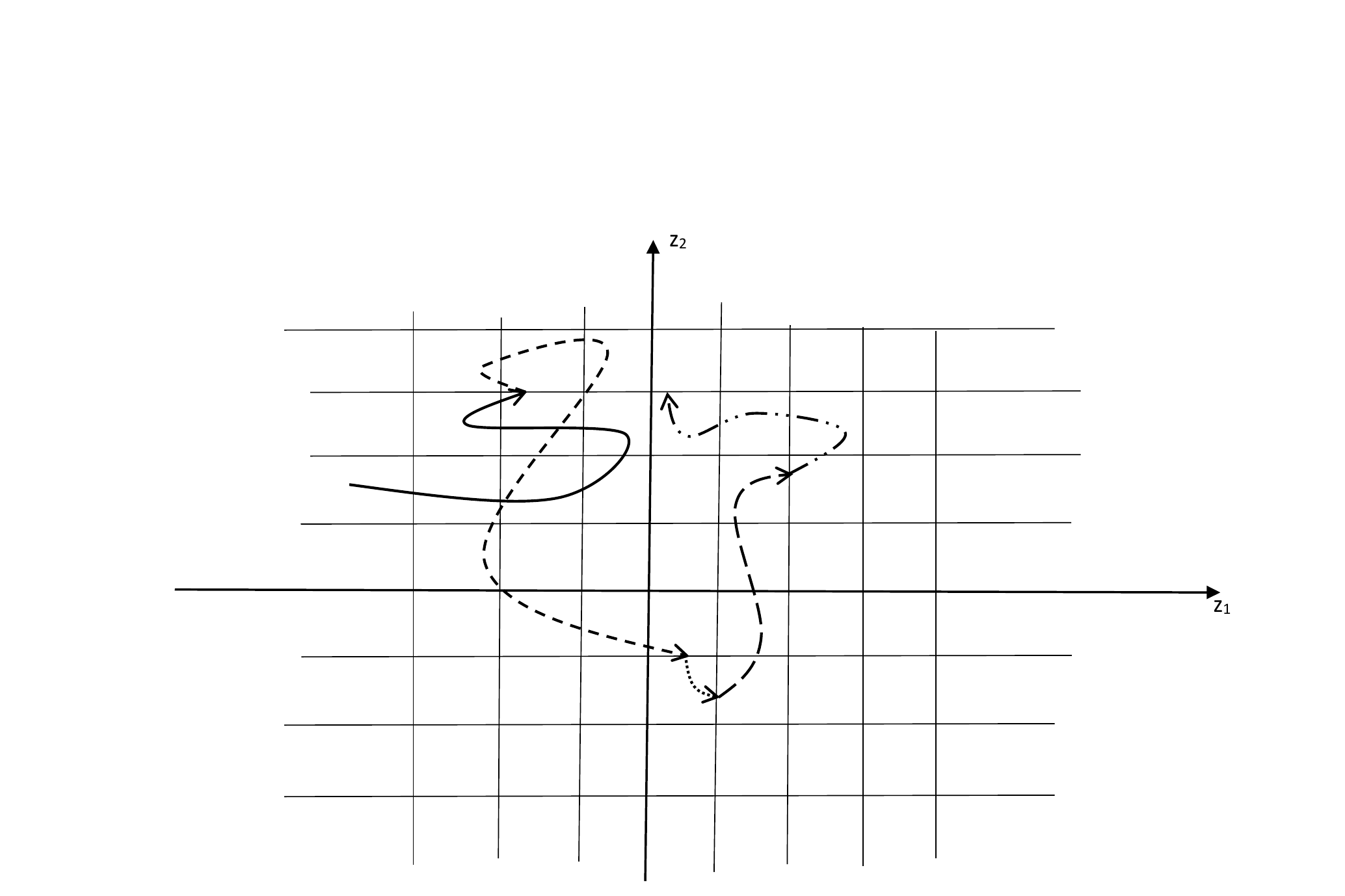}
}
\caption{\label{fig:multiswitching} Multiswithcing dynamics in the case of $m=2$ fields, and $k=4$ relays (compare with Figure \ref{fig:play_discreto}), and with $\zeta_1=e_1$, $\zeta_2=e_2$ the first and second vectors of the canonical basis of $\mathbb{R}^n$. Here (compare with Figure \ref{projection})  it is represented the projection on the plane $(z_1,z_2)$ of a possible trajectory starting from $z_1<-3/4$, $1/4<z_2<3/4$ with $w_4[z_1](0)=-1$ corresponding to $(-1,-1,-1,-1)$, and $w_4[z_2](0)=0$ corresponding to $(-1,-1,1,1)$. The filled trajectory evolves by $(g_1^{(-1,-1,-1,-1)},g_2^{(-1,-1,1,1)})$. The dashed trajectory evolves by $(g_1^{(-1,-1,-1,-1)},g_2^{(-1,1,1,1)})$. The pointed trajectory evolves by $(g_1^{(-1,-1,-1,-1)},g_2^{(-1,-1,1,1)})$. The long-dashed trajectory evolves by $(g_1^{(-1,-1,-1,1)},g_2^{(-1,-1,1,1)})$. The dashed-double-pointed trajectory evolves by $(g_1^{(-1,-1,1,1)},g_2^{(-1,-1,1,1)})$.}
\end{figure}

{\it Proof.}
We can restrict to the case with an initial output states satisfying (\ref{eq:hysteron_k}). Then for any $j=1,\dots,m$ the evolution of the pair $(z\cdot\xi_j,w_k[z\cdot\xi_j])$ is inside the discrete hysteresis-loop of Figure \ref{fig:play_discreto} and $w_k[z\cdot\xi_j]$ is uniquely generated by the admissible $k$-tuple $(h_1[z\cdot\xi_j],\dots,h_k[z\cdot\xi_j])=(1,\dots,1,-1,\dots,-1)$. Hence, for any $j$ and any one of those $k$-tuple (which, by (\ref{eq:hysteron_k}) and our choice of the thresholds, are exactly $k+1$: $(-1,-1,\dots,-1),(1,-1,\dots,-1),\dots,(1,\dots,1)$) we have a field (see Figure \ref{fig:multiswitching})

\[
z\mapsto g_j^{(h_1,\dots,h_k)}(z)
\]

The proof can be made as in Theorem \ref{thm:switching}. Indeed, also in the case of superposition of $k$ switchings, the strategy used to reach the final point depends on where the initial one is located. More precisely, even in this multiple switching situation, we can distinguish different intersecting sectors, no-one with empty interior, so that any point in $\mathbb{R}^n$ stays in the interior of a sector, (see example in Figure \ref{fig:multiswitching}). Thus if the initial an final points are in the interior of the same sector we can use the strategy of case $1$ of Theorem \ref{thm:switching}, if instead one end-point is in the interior of a sector and the other end-point is in the interior of another sector, then the strategy can be the one of case 2 of Theorem \ref{thm:switching}, maybe switching more than once. Finally if the starting point is on the boundary of a sector the strategy will be the same of case $3$ of Theorem \ref{thm:switching}. All the strategies make use of controllability Assumption \eqref{eq:hypothesis_versus_truncated}.
 \EndProof

\begin{remark}
For a result of controllability of (\ref{eq:play_truncated_system}), one should pass to the limit in $k\to+\infty$, in the controllabilty problem (\ref{eq:discrete_system}). This will be the subject of future studies. 
Here we note that the problem (\ref{eq:play_truncated_system}) seems to be not exactly fitting the similar problem in (\ref{eq:intro})-right, because of the explicitly presence of the variable $z$ inside the fields with hysteresis. Moreover, in this case, a triangular feature as in Subsections 41. and 4.2 seems to be not necessary. One crucial point is the controllabilty hypotheses (\ref{eq:hypothesis_versus_truncated}). 
That hypothesis can be rather natural in some cases. Let us consider the system (without hysteresis)

\begin{equation}
\label{eq:1}
\dot z=\sum_{j=1}^m\tilde g_j(z)u_j
\end{equation}

\noindent
where the fields $\tilde g_j$ satisfy the controllability conditions (Chow). Actually, due for example to some kind of damage, we do not exactly face that system, but a perturbation of it of the form

\[
\dot z=\sum_{j=1}^m\left(\tilde g_j(z)+f_j({\cal P}[z\cdot\zeta_j])\right)u_j
\]

\noindent
where $f_j:\mathbb{R}\to\mathbb{R}^n$. This can be seen as a generalization of a linear system with feedback control which is affected by some damage (see Tarbouriech et al. \cite{Tarbouriech} and Visintin \cite{Visintin94} for more details on damaged systems and hysteresis.)
 Another possible model is

\begin{equation}
\label{eq:2}
\dot z=\sum_{j=1}^mf_j({\cal P}[z\cdot\zeta_j])\tilde g_j(z)u_j
\end{equation}

\noindent
where $f_j:\mathbb{R}\to\mathbb{R}$. More generally

\begin{equation}
\label{eq:3}
\dot z=\sum_{j=1}^m g_j(z,{\cal P}[z\cdot\zeta_j])u_j
\end{equation}

\noindent
where $g_j:\mathbb{R}^{n+1}\to\mathbb{R}^n$.
Again the perturbation can be seen as a damage in the feedback control.
Since systems (\ref{eq:1}) is controllable, under some reasonable hypotheses on the perturbation $f_j$,
the presence of the hysteretic term does not affect the controllability of the "non-perturbed" part, that is (\ref{eq:hypothesis_versus_truncated}) holds.

When we perform the discrete approximation of the play operator (\ref{eq:w_n}), as already explained, if we are in the hysteresis-loop as in Figure \ref{fig:play_discreto}, to any possible value $w_k^j$ of the output of $w_k[z\cdot\zeta_j]$ a unique suitable $k$-string, $s^j_k\in S_k$, of $1$ and $-1$ is associated (where $S_k$ is the set of all such suitable $k$-string).
In any sector where the string does not change (see Figure \ref{fig:multiswitching}), we then move with the fields $z\mapsto g_j^{s^j_k}(z)=g_j(z,w_j^k)$. By (\ref{eq:hypothesis_versus_truncated}), for every choice of the strings $s^j_k$, the fields $z\mapsto g_j^{s_j^k}(z)$ satisfy the controllability condition.

Note that, for example in the case (\ref{eq:1}), for every fixed string $s_j^k$, it is $\nabla g_j^{s_j^k}(z)=\nabla g_j(z)$. So (\ref{eq:hypothesis_versus_truncated}), which involves the Lie brackets of the vector fields and thus their derivatives, is not so unrealistic in the case, for example, of small magnitude perturbation $f_j$. Similar considerations may be done in the case (\ref{eq:2}). 
Moreover also note that (\ref{eq:hypothesis_versus_truncated}) does not take care of the evolution of the perturbation variable $t\mapsto w_j(t)={\cal P}[z\cdot\zeta_j](t)$, which certainly may further affect the controllability of (\ref{eq:3}). On the contrary, the discrete problem makes use of (\ref{eq:hypothesis_versus_truncated}) only, because it leads to use it in the sectors where the variable $w^j_k$ (as well as the $k$-string $s_k^j$) does not evolve, but, on the other side, when $k$ increases, the number of those regions also increases and they also present finer granularity. The passage to the limit $k\to+\infty$ it is certainly worth studying.
\end{remark}


\end{document}